\documentclass[notitlepage,11pt]{article}
\usepackage{amssymb,amsmath,geometry%, cases
}

\catcode`\@=11
\@addtoreset{equation}{section}

\catcode`\@=12

\allowdisplaybreaks

\usepackage{color}

\def\R{\mathbb{R}}

\def\f{\varphi}
\def\iirn{\iint\limits_{\R^n\!\times\R^n}}
\def\iirnp{\iint\limits_{\R^n_+\!\times\R^n_+}}

\def\irn{\int\limits_{\R^n}}

\def\irnplus{\int\limits_{\R^n_+}}
\def\irnminus{\int\limits_{\R^n_-}}

\def\eps{\varepsilon}

\def\Dplus{-\Delta_{\R^n_+}^{\!N}\!}  %Standard s=1 half space

\def\Ds{\left(-\Delta\right)^{\!s}\!}  %Stanfard Dirichlet-Fourier
\def\DsNROmega{(-\Delta_{\Omega}^{\!N})^{\!s}_{\!{\rm R}}}  %Neumann Restricted
\def\DsNR{(-\Delta_{\R^n_+}^{\!N}\!)^{\!s}_{\!{\rm R}}}  

\def\DsNSOmega{(-\Delta_{\Omega}^{\!N})^{\!s}_{\!{\rm Sr}}}  %Neumann Semirestricted
\def\DsNS{(-\Delta_{\R^n_+}^{\!N}\!)^{\!s}_{\!{\rm Sr}}}  

\def\DsNSpOmega{(-\Delta_{\Omega}^{\!N})^{\!s}_{\!{\rm Sp}}}  %Neumann Spectral
\def\DsNSp{(-\Delta_{\R^n_+}^{\!N}\!)^{\!s}_{\!{\rm Sp}}}

\def\sstar{{2^*_s}}
\def\sigstar{{2^*_\sigma}}

\def\proof{\noindent{\textbf{Proof. }}}
\def\QED{\hfill {$\square$}\goodbreak \medskip}

\newtheorem{Theorem}{Theorem}[section]
\newtheorem{Lemma}[Theorem]{Lemma}

\newtheorem{Corollary}[Theorem]{Corollary}
\newtheorem{Remark}[Theorem]{Remark}

\begin{document}

\title %{\vspace{-10mm}
{Sobolev inequalities for Neumann Laplacians on half spaces}

\author{Roberta Musina\footnote{Dipartimento di Matematica ed Informatica, Universit\`a di Udine,
via delle Scienze, 206 -- 33100 Udine, Italy. Email: {roberta.musina@uniud.it}. 
{Partially supported by Miur-PRIN project
''Variational methods, with applications to problems in mathematical physics and geometry'' 2015KB9WPT\_001.}}~ and
Alexander I. Nazarov\footnote{
St.Petersburg Department of Steklov Institute, Fontanka 27, St.Petersburg, 191023, Russia, 
and St.Petersburg State University, 
Universitetskii pr. 28, St.Petersburg, 198504, Russia. E-mail: al.il.nazarov@gmail.com.
Partially supported by  RFBR grant 17-01-00678. 
%RFBR grant 14-01-00534 and by St.Petersburg University grant 6.38.670.2013.
}
}

\date{}

\maketitle

\begin{abstract}
{\small 
We consider different fractional Neumann Laplacians %$(-\Delta^N_\Omega)^s$ 
of order $s\in(0,1)$,  namely,
the {\em Restricted Neumann Laplacian} $\DsNROmega$, 
 the {\em Semirestricted Neumann
Laplacian} $\DsNSOmega$ and the {\em Spectral Neumann Laplacian} $\DsNSpOmega$.
In particular, we are interested in attainability of Sobolev constants for these operators in $\Omega=\R^n_+$.}

\medskip

\medskip

\noindent
\textbf{2010 Mathematics Subject Classfication:} 47A63; 35A23.
\end{abstract}

\normalsize

\bigskip

\section{Introduction}
\label{S:Introduction}
Let $n\ge 3$ be an integer and put $2^*=\frac{2n}{n-2}$. 
By the Sobolev inequality, the Hilbert  space 
$$
{\mathcal D}^1(\R^n)=\{u\in {L^{2^*}(\R^n)}~|~\langle-\Delta u,u\rangle=\displaystyle\irn|\nabla u|^2~\!dx<\infty~\},
$$
is continuously embedded into $L^{2^*}(\R^n)$. 
It has been proved in \cite{Au, Ta} that
the radial function $U(x)=(1+|x|^2)^{\frac{2-n}{2}}$ achieves the Sobolev constant 
\begin{equation}
\label{eq:Sob1}
\mathcal S=\inf_{u\in {\mathcal D}^1(\R^n)\atop\scriptstyle  u\ne 0}
\frac{\langle-\Delta u,u\rangle}{\|u\|^2_{L^{2^*}(\R^n)}}
\end{equation}
(for $n=3$ this remarkable fact was established earlier in \cite{Rosen}).

Next, let $\R^n_+$ be an half-space, for instance 
$$\R^n_+=\{~\!x=(x_1,x')\in\R\times\R^{n-1}~|~~ x_1>0~\!\}.$$
We denote by $\Dplus$ the Neumann Laplacian in $\R^n_+$, that is the distribution
$$\displaystyle{\langle \Dplus u,\f\rangle=\irnplus\nabla u\cdot\nabla\f~\!dx}$$
for $u,\f\in \mathcal D^1(\R^n_+)=\{u\in {L^{2^*}(\R^n_+)}~|~\langle \Dplus u,u\rangle=\irnplus|\nabla u|^2~\!dx<\infty~\}.$
If $u$ is smooth enough we clearly have $\Dplus u=-\Delta u$ pointwise on $\R^n_+$.
It is easy to check that the Aubin-Talenti function $U$ solves the Neumann problem
$$
-\Delta U= n(n-2) U^{2^*-1}\quad\text{in $\R^n_+$~,}\quad \frac{\partial U}{\partial x_1}=0\quad\text{on $\partial \R^n_+$},
$$
and achieves the Neumann Sobolev constant
$$
\mathcal S(\R^n_+)=\inf_{u\in {\mathcal D}^1(\R^n_+)\atop\scriptstyle  u\ne 0}
\frac{\langle \Dplus u,u\rangle}{\|u\|^2_{L^{2^*}(\R^n_+)}}~\!.
$$
In particular, one infers that 
$\mathcal S(\R^n_+)=2^{-\frac2n}\mathcal S$. This crucial observation permits to relate 
existence/multiplicity phenomena in 
critical/nearly critical Neumann problems on %domains 
$\Omega\subset\R^n$
to the geometric properties of $\partial\Omega$. 
Due to the abundant literature on this
subject, we limit ourselves to cite the the
pioneering results in  \cite{AM1,  AMY}, the more recent papers 
\cite{delMM, WWY}, the surveys \cite{KN, Naz-surv} and references therein.

\medskip

The goal of the present paper is to study Sobolev-type constants on half spaces governed by
Neumann fractional Laplacians of order
$s\in(0,1)$.

\bigskip

We will discuss three different nonlocal operators, namely,
the {\em Restricted Laplacian} $\DsNR$ (in Section \ref{S:Restricted}), the
{\em Spectral Laplacian} $\DsNSp$ (in Section \ref{S:Spectral}), and the {\em Semirestricted Laplacian} $\DsNS $
(in Section \ref{S:Semirestricted}). 

We always assume  $s\in (0,1)$ and
$n>2s$, that  is  a restriction only if $n=1$, and  put
$$
\sstar:=\frac{2n}{n-2s}~\!.
$$
Before describing our main results we
recall some facts about the {\em Dirichlet Laplacian}
$$
\Ds u(x)={C_{n,s}}~\!\cdot~\!{\rm P.\!V.}\irn\frac{u(x)-u(y)}{|x-y|^{n+2s}}~dy~,%\quad x\in\R^n~\!,
\quad C_{n,s}=\displaystyle{\frac{s2^{2s}\Gamma\big(\frac{n}{2}+s\big)}{\pi^{\frac{n}{2}}\Gamma\big(1-s\big)}}~\!.
$$
Here $u$ is a smooth and rapidly decreasing function on $\R^n$, P.\!V. means principal value and $x$ runs in the whole space $\R^n$.

Thanks to the Hardy-Littlewood-Sobolev inequality, the quadratic form $\langle \Ds u,u\rangle$ 
induces an Hilbertian structure on the space
$$
{\mathcal D}^s(\R^n)\!=\!\{u\in {L^\sstar(\R^n)}~|~\langle \Ds u,u\rangle=\frac{C_{n,s}}{2}\!\iirn\frac{(u(x)-u(y))^2}{|x-y|^{n+2s}}~dxdy<\infty~\},
$$
that is continuoulsy embedded into ${\mathcal D}^s(\R^n) \hookrightarrow L^\sstar(\R^n)$.
It has been proved in \cite{CoTa}
that
\begin{equation}
\label{eq:CoTa}
\mathcal S_s:=\inf_{\scriptstyle u\in \mathcal D^s(\R^n)\atop\scriptstyle  u\ne 0}
\frac{\langle \Ds u,u\rangle}
{\|u\|^2_{L^\sstar(\R^n)}}=
2^{2s}\pi^{2s}~\frac{\Gamma\big(\frac{n}{2}+s\big)}
{\Gamma\big(\frac{n}{2}-s\big)}~\!
\left[\frac{\Gamma\big(\frac{n}{2}\big)}{\Gamma({n})}\right]^{\frac{2s}{n}}
\end{equation}
and that, up to dilations, translations and multiplications, the fractional Sobolev constant
$\mathcal S_s$ is attained only by the function
\begin{equation}\label{U}
U_{\!s}(x)=\left(1+|x|^2\right)^{\frac{2s-n}{2}}.
\end{equation}
We are now in position to describe the Neumann Laplacians we are interested in.

\medskip

The {\em Restricted (or {Regional}) fractional Laplacian} on the half space $\R^n_+$
is formally defined by
$$
\DsNR u(x)={C_{n,s}}~\!\cdot~\!{\rm P.\!V.}\irnplus\frac{u(x)-u(y)}{|x-y|^{n+2s}}~dy~,\quad x\in\R^n_+~\!.
$$
Restricted Laplacians appear as generators of so-called censored processes. 
A large number of papers deal with operators $\DsNROmega$ on domains $\Omega\subset \R^n$;
we limit ourselves to cite  \cite{DZ, G1, G2, W1, W2, W3} and references therein.

In Lemma \ref{L:spaceR} we prove that the {\em Restricted quadratic form} $\langle \DsNR u,u \rangle$
induces a Hilbertian norm on the space
$$
{\mathcal D}^s_{\rm R}(\R^n_+)=\left\{
u\in L^\sstar(\R^n_+)~|~  \langle \DsNR u,u \rangle= 
\frac{C_{n,s}}{2}\iirnp \frac{(u(x)-u(y))^2}{|x-y|^{n+2s}}~\!\!dxdy<\infty~\right\},
$$
and that ${\mathcal D}^s_{\rm R}(\R^n_+)$ is continuously embedded into $L^\sstar(\R^n_+)$. Hence, 
the  {\em Restricted Sobolev constant}
\begin{equation}
\label{eq:problemR}
S^{\rm R}_s(\R^n_+)=\inf_{u\in {\mathcal D}^s_{\rm R}(\R^n_+)\atop\scriptstyle  u\ne 0}
\frac{\langle \DsNR u,u \rangle}{\|u\|^2_{L^\sstar(\R^n_+)}}
\end{equation}
is positive. In Section \ref{S:Restricted} we prove the following existence theorem. 

\begin{Theorem}
\label{T:Restricted}
 It holds that $\mathcal S^{\rm R}_s(\R^n_+)<2^{-\frac{2s}{n}}\mathcal S_s$, and $\mathcal S^{\rm R}_s(\R^n_+)$
is achieved. 
\end{Theorem}

The {\em Neumann Spectral  fractional Laplacian} $\DsNSp$ is the $s$-th power of the standard Neumann Laplacian in the sense of spectral theory. 
In Section \ref{S:Spectral} we prove the next existence result.

\begin{Theorem}
\label{T:Spectral}
It holds that 
$$
\mathcal S^{\rm Sp}_s(\R^n_+):=\inf_{u\in {\mathcal D}^s_{\rm R}(\R^n_+)\atop\scriptstyle  u\ne 0}
\frac{\langle\DsNSp u,u\rangle}{\|u\|^2_{L^\sstar(\R^n_+)}}
=2^{-\frac{2s}{n}}\mathcal S_s~\!,
$$
and $\mathcal S^{\rm Sp}_s(\R^n_+)$
is achieved by the function $U_s\big|_{\R^n_+}$. 
\end{Theorem}

The {\em Semirestricted Laplacian}  $\DsNSOmega$ is the operator
$$
\DsNS u=\chi_{\R^n_+}(-\Delta)^{\!s} u+\chi_{\R^n_-}\DsNR u ~,\quad x\in\R^n~\!.
$$
For general domains $\Omega\subset\R^n$, the Semirestricted Laplacian $\DsNSOmega$
can be used to study non-homogeneous Dirichlet problems for $\Ds$ on $\Omega$,
see for instance the survey paper \cite{Ro} by Ros-Oton, and has been proposed by  Dipierro, Ros-Oton and Valdinoci
\cite{DRoV} as an alternative approach to Neumann problems. 

By the computations in \cite[Lemma 3]{DRoV}, one naturally associates to 
$\DsNS$ the {\em Semirestricted quadratic form}
$$
\langle\DsNS u,u\rangle=\frac{C_{n,s}}{2} \iint\limits_{\R^{2n}\setminus(\R^n_-)^2} \frac{(u(x)-u(y))^2}{|x-y|^{n+2s}}~\!\!dxdy~\!.
$$
In Lemma \ref{L:Sr_space} we show that 
$\langle\DsNS u,u\rangle$ induces an Hilbertian structure on the space
 \begin{equation}
\label{eq:SR_space}
{\mathcal D}^s_{\rm Sr}(\R^n_+)=\{~\!\text{measurable $u:\R^n\to \R$}~\!\big|~\!
u|_{\R^n_+}\in L^\sstar(\R^n_+),~\langle\DsNS u,u\rangle<\infty~\!\}
\end{equation}
and that ${\mathcal D}^s_{\rm Sr}(\R^n_+)$ is continuously embedded into $L^\sstar(\R^n_+)$. Thus the
{\em Semirestricted Sobolev constant}
\begin{equation}
\label{eq:DRoV}
\mathcal S^{\rm Sr}_s(\R^n_+):=\inf_{u\in {\mathcal D}^s_{\rm Sr}(\R^n_+)\atop u|_{\R^n_+}\neq 0}%{\mathcal D}^s(\R^n)}
\frac{ \langle\DsNS u,u\rangle}{\|u\|^2_{L^\sstar(\R^n_+)}}
\end{equation}
is positive. 
In Section \ref{S:Semirestricted}  we prove the next theorem.

\begin{Theorem}
\label{T:Semirestricted}
It holds that $\mathcal S^{\rm Sr}_s(\R^n_+)\le \mathcal S_s$. If
$\mathcal S^{\rm Sr}_s(\R^n_+)<\mathcal S_s$ then $\mathcal S^{\rm Sr}_s(\R^n_+)$ is achieved.
\end{Theorem}

Some sufficient conditions for the
validity of the strict inequality $\mathcal S^{\rm Sr}_s(\R^n_+)<\mathcal S_s$ are proved in Theorem \ref{T:strictSr}.

In Section \ref{S:subcr} we collect some results on attainability of sharp constants in Hardy-Sobolev inequalities for all considered Neumann 
Laplacians.

The proofs of our main theorems need few results, some of which are of independent interest.
In Appendix A we prove some auxiliary properties of the ``best extension'' projector $P_{\!s}$ arising when studying the Semirestricted 
Laplacian. In Appendix B we study the limit properties of quadratic forms as $s\to0$ and as $s\to1$.

\bigskip
\noindent
{\small 
{\bf Notation.}
For $Z\subseteq\R^n\times\R^n$ and for any function $u$ such that $(x,y)\mapsto u(x)-u(y)$ is a measurable function on $Z$  we put 
$$
\mathcal E_s(u;Z)=\frac{C_{n,s}}{2}\iint\limits_Z\frac{(u(x)-u(y))^2}{|x-y|^{n+2s}}~dxdy.
$$
In particular, we  have $\langle \Ds u,u\rangle=\mathcal E_s(u;\R^n\!\times\R^n\!)$ and
$$
\langle \DsNR u,u\rangle=\mathcal E_s(u;\R^n_+\!\times\R^n_+\!)~,\quad
\langle \DsNS u,u\rangle=\mathcal E_s(u;\R^{2n}\setminus(\R^n_-)^2).
$$
We set
$$
\R^n_{\pm}=\{x\in\R^n\,\big|\, x_1\gtrless0\}; \qquad B_\rho(y)=\{x\in\R^n\,\big|\,|x-y|<\rho\}; \qquad B_\rho=B_\rho(0).
$$
For a domain $\Omega\subset\R^n$ we put $\Omega^{\pm}=\Omega\cap\R^n_{\pm}$.

Through the paper, all constants depending only on $n$ and $s$ are denoted by $c$. To indicate that a constant depends on other quantities %$a,b,\dots$ 
we list them in parentheses: $c(\dots)$.}

\section{Restricted Laplacian and proof of Theorem \ref{T:Restricted}}
\label{S:Restricted}
A few preliminaries are in order.

First, notice that ${\mathcal D}^s(\R^n)\cap L^2(\R^n)$ is the standard Sobolev space
$H^s(\R^n)$ (we refer to \cite{Tr} for basic results about $H^s$-spaces). In particular
${\mathcal D}^s(\R^n)\supsetneq H^s(\R^n)$ and ${\mathcal D}^s(\R^n)$ is a subset of $H^s_{\rm loc}(\R^n)$,
that means $\f u\in H^s(\R^n)$ for  $\f\in\mathcal C^\infty_0(\R^n)$ and $u\in {\mathcal D}^s(\R^n)$.
Therefore, 
$\mathcal C^\infty_0(\R^n)$ is dense in ${\mathcal D}^s(\R^n)$ and
${\mathcal D}^s(\R^n)$ is compactly embedded into 
$L^p_{\rm loc}(\R^n)$ for any $p\in[1,\sstar)$.

For future convenience we introduce $ \widetilde{\mathcal D}^s(\R^n_+)$ as the closure in 
${\mathcal D}^s(\R^n)$ of $C^\infty_0(\R^n_+)$. By standard arguments and direct computation one can check that 
\begin{gather}
\nonumber
 \widetilde{\mathcal D}^s(\R^n_+)
=\big\{u\in \mathcal D^s(\R^n)~|~u\equiv 0~~\text{on $\R^n_-$}~\big\},\\
\label{eq:N_sum}
\mathcal E_s(u;\R^n\!\!\times\!\R^n)=\mathcal E_s(u;\R^n_+\!\times\!\R^n_+)+\gamma_s\irnplus x_1^{-2s}u^2~\!dx
\quad\text{for any $u\in \widetilde{\mathcal D}^s(\R^n_+)$,}
\end{gather}
where
\begin{equation}
\label{eq:gamma}
\gamma_s:=C_{n,s}|x_1|^{2s}\irnplus\frac{dy}{|x-y|^{n+2s}}=\frac{2^{2s-1}\Gamma\Big(s+\frac12\Big)}{\sqrt\pi~\!\Gamma\big(1-s\big)}
\quad\text{for any $x\in\R^n_-$.}
\end{equation}
We are now in position to start our description of the space $\mathcal D^s_{\rm \!R}(\R^n_+)$.

\begin{Lemma}
\label{L:spaceR}
\begin{itemize}
\item[$i)$] 
$\|u\|^2_{\mathcal D^s_R(\R^n_+)}:=\mathcal E_s(u;\R^n_+\!\times\!\R^n_+)$
is an Hilbertian norm on $\mathcal D^s_{\rm \!R}(\R^n_+)$ and for $u\in \mathcal D^s_{\rm R}(\R^n_+)$ the distribution
$\DsNR u$ defined by 
$$
\langle \DsNR u,v\rangle= \frac{C_{n,s}}{2}\iirnp \frac{(u(x)-u(y))(v(x)-v(y))}{|x-y|^{n+2s}}~\!\!dxdy~,\quad
v\in \mathcal D^s_{\rm R}(\R^n_+)
$$
belongs to the dual space $\mathcal D^s_{\rm R}(\R^n_+)'$;
\item[$ii)$] $\mathcal D^s_{\rm R}(\R^n_+)$ is continuously embedded into $L^\sstar(\R^n_+)$;
\item[$iii)$] $\mathcal D^s_{\rm R}(\R^n_+)$ is compactly embedded in $L^p_{\rm loc}(\overline{R}^n_+)$, for any $p\in[1,\sstar)$. 
That is,  any bounded sequence in $\mathcal D^s_{\rm R}(\R^n_+)$ has a subsequence that converges in $L^p(\Omega^+)$, for any
bounded domain $\Omega\subset\R^n$;
\item[$iv)$]
$\widetilde{\mathcal D}^s(\R^n_+)$ is continuously embedded into $\mathcal D^s_{\rm R}(\R^n_+)$.
\end{itemize}
\end{Lemma}

\proof
For any $u\in \mathcal D^s_{\rm R}(\R^n_+)$ we denote by
$\hat u:\R^n\to \R$ the even extension of $u$,
that is,
\begin{equation}
\label{eq:even_ex}
\hat u(x_1,x')=u(|x_1|,x')~\!.
\end{equation}
Clearly $\hat u\in L^\sstar(\R^n)$ and $\|\hat u\|^\sstar_{L^\sstar(\R^n)}=2
\| u\|^\sstar_{L^\sstar(\R^n_+)}$. Moreover, 
we have
\begin{equation}
\label{eq:numero0}
\begin{cases}
\mathcal E_s(\hat u;\R^n\!\!\times\!\R^n)=2\mathcal E_s(u;\R^n_+\!\times\!\R_+^n)+2\mathcal E_s(\hat u;\R^n_+\times\R_-^n)\\
\mathcal E_s(\hat u;\R^n_+\times\R_-^n)\le \mathcal E_s(u;\R^n_+\!\times\!\R_+^n)~\!.
\end{cases}
\end{equation}
Hence $\hat u\in \mathcal D^s(\R^n)$ and using also the Sobolev inequality we get
$$
\mathcal E_s(u;\R^n_+\!\times\!\R_+^n)\ge \frac14 \mathcal E_s(\hat u;\R^n\!\!\times\!\R^n)\ge 
\frac14\mathcal S_s \|u\|^2_{L^\sstar(\R^n_+)}\quad\text{for any $u\in \mathcal D^s_{\rm R}(\R^n_+)$.}
$$
We infer that $\mathcal E_s(~\!\cdot~\!;\R^n_+\!\times\!\R_+^n)^{\frac12}$ is a norm on $\mathcal D^s_{\rm R}(\R^n_+)$ and that
$\mathcal D^s_{\rm R}(\R^n_+)$ is continuously embedded into $ L^\sstar(\R^n_+)$. 

The conclusion of the proof easily follows
from the continuity of the operators $u\mapsto \hat u$, $\mathcal D^s_{\rm R}(\R^n_+)\to \mathcal D^s(\R^n)$ and 
$u\mapsto u|_{\R^n_+}$, $\mathcal D^s(\R^n)\to \mathcal D_{\rm R}^s(\R^n_+)$.
\QED

\begin{Lemma}
\label{L:technicalR}
Let $u\in \mathcal D^s_{\rm R}(\R^n_+)$, $\f\in\mathcal C^\infty_0(\R^n)$ and let $\Omega
\subset\R^n$ be a bounded
 domain containing the support of $\f$. Then $\f u\in \mathcal D^s_{\rm R}(\R^n_+)$
and 
$$
\big|\mathcal E_s(\f u;\R^n_+\!\!\times\!\R^n_+)-\langle \DsNR u,\f^2 u\rangle\big|
\le  c(\f) \|u\|_{L^2(\Omega^+)}\left(
\|u\|_{L^2(\Omega^+)}
+\mathcal E_s(u;\R^n_+\!\times\R^n_+)^\frac12\right).
$$
\end{Lemma}

\proof
This is an adaptation of \cite[Lemma 2.1]{MNcn}; we restrict ourselves to indicate the main changes in the proof.

To simplify notation we put
\begin{equation}
\label{eq:to_simplify}
\Psi_\f(x,y)=\frac{(\f(x)-\f(y))^2}{|x-y|^{n+2s}}~.
\end{equation}
As in \cite[Lemma 2.1]{MNcn} we estimate 
\begin{equation}
\label{eq:Psi_estimate}
\irn\Psi_{\!\f}(x,y)~\!dy\le c(\f),
\end{equation}
with $c(\f)$ not depending on $x$, that readily gives $\f u\in \mathcal D^s_{\rm R}(\R^n_+)$ via standard arguments.
Next, by direct computation one finds 
$$
\mathcal E_s(\f u;\R^n_+\!\!\times\!\R^n_+)-\langle\DsNR u,\f^2 u\rangle=
c\iirnp{u}(x){u}(y)\Psi_\f(x,y)~\!dxdy=:\mathcal B_\f~\!.
$$
Since the support of $\Psi_\f$ is contained in $(\Omega\times\R^n)\cup(\R^n\times\Omega)$, we have
$$
c~\!\big|\mathcal B_\f\big|\le \iint\limits_{\Omega^+\!\times\Omega^+} |u(x)u(y)|\Psi_\f(x,y)~\!dxdy+
\int\limits_{\Omega^+}|u(x)|\Big(\int\limits_{\R^n_+\setminus\Omega} \frac{|u(y)||\f(x)|^2}{|x-y|^{n+2s}}dy
\Big)dx~\!.
$$
We use the triangle inequality $|u(y)|\le |u(x)|+ |u(x)-u(y)|$ in the last integral to infer
$c\big|\mathcal B_\f\big|\le  I_1+I_2+\|\f\|_\infty I_3$, where
\begin{gather*}
I_1=\iint\limits_{\Omega^+\times\Omega^+} |u(x)u(y)|\Psi_\f(x,y)~\!dxdy~,\quad
I_2=\int\limits_{\Omega^+} |u(x)\f(x)|^2\Big(\int\limits_{\R^n_+\setminus\Omega}\frac{dy}{|x-y|^{n+2s}}\Big)dx~,\\
I_3=\iint\limits_{\Omega^+\times (\R^n_+\setminus\Omega)} 
\frac{|u(x)-u(y)|}{|x-y|^{\frac{n+2s}{2}}}~\!\frac{|u(x)\f(x)|}{|x-y|^{\frac{n+2s}{2}}}~\!dxdy.
\end{gather*}
Arguing  as in the proof of Lemma 2.1 of \cite{MNcn} (with $\Omega$ replaced by $\Omega^+$)
one gets the estimates
\begin{gather*}
I_1\le  c(\f)\int\limits_{\Omega^+}|u(x)|^2~\!dx~,\qquad
I_2\le c(\f) \int\limits_{\Omega^+}|u(x)|^2dx,\\
I_3^2\le c ~\! \mathcal E_s(u;\R^n_+\!\times\R^n_+) I_2 
\le c(\f) \mathcal E_s(u;\R^n_+\!\times\R^n_+) \int\limits_{\Omega^+}|u(x)|^2dx~\!,
\end{gather*}
that end the proof. 
\QED

\bigskip
\noindent 
{\bf Proof of Theorem \ref{T:Restricted}.}
We test $\mathcal S^{\rm R}_s(\R^n_+)$ with the function ${U_{\!s}}|_{\R^n_+}$, see (\ref{U}). Since
$U_{\!s}$ is even in $x_1$ we have 
$\frac12\mathcal E_s(U_s;\R^n\!\!\times\!\R^n)=\mathcal E_s(U_{\!s};\R^n_+\!\times\!\R_+^n)+\mathcal E_s(U_{\!s};\R^n_+\times\R_-^n)
>\mathcal E_s(U_{\!s};\R^n_+\!\times\!\R_+^n)$.
Thus
$$
\mathcal S_s=\frac{\mathcal E_s(U_{\!s};\R^n\!\times\!\R^n)}{\|U_{\!s}\|^2_{L^\sstar(\R^n)}}>
2^{\frac{2s}{n}}~\!\frac{\mathcal E_s(U_{\!s};\R^n_+\!\times\!\R^n_+)}{\|U_{\!s}\|^2_{L^\sstar(\R^n_+)}}
\ge 2^{\frac{2s}{n}}S^{\rm R}_{s}(\R^n_+),
$$
and the first claim is proved. 

Now we show that the noncompact minimization problem (\ref{eq:problemR}) 
admits a solution. 
We follow the outline of the proof of Theorem 0.1 in \cite{GM}, see also \cite{MNcn}.

Thanks to  a standard convexity argument, we only need to construct a bounded minimizing sequence
$u_h\in \mathcal D^s_{\rm R}(\R^n_+)$ such that $u_h\to u\neq 0$ weakly. 
We put 
$${\mathcal S_R}=\mathcal S^{\rm R}_s(\R^n_+)~,\qquad E(u)=\mathcal E_s(u;\R^n_+\!\!\times\!\R^n_+)$$
and we limit ourself to the more difficult case $n\ge 2$. For $\rho>0$ and ${z}\in \R^{n-1}$ we denote by $B'_\rho({z})$ the  $(n-1)$-dimensional ball
$$
B'_\rho({z})=\{x'\in\R^{n-1}~|~|x'-{z}|<\rho~\}~\!.
$$
Then we take a small number $\eps_0$ and a finite number of points  
$x'_1,\cdots x'_\tau\in \R^{n-1}$ such that
\begin{equation}
\label{eq:Rcover}
0<\eps_0<\frac12 {\mathcal S_R}~,\quad \overline{B'_2(0)}\subset \bigcup_{j=1}^\tau B_{1}'(x'_j).
\end{equation}
Since the ratio in (\ref{eq:problemR}) is invariant with respect to translations in $\R^{n-1}$ and
with respect to the transforms ${\mathcal D}^s_{\rm R}(\R^n_+)\to
{\mathcal D}^s_{\rm R}(\R^n_+)$,
$u(x)\mapsto  \alpha u(\beta x)$ ( $\alpha\neq 0, \beta>0$), we can find
a bounded minimizing sequence $u_h$ for $\mathcal S_R$ such that
\begin{gather}
\label{eq:normaliz}
\|u_h\|^\sstar_{L^\sstar(\R^n_+)}={\mathcal S_R^\frac{n}{2s}}~,~~
{E}(u_h)={\mathcal S_R^\frac{n}{2s}}+o(1)
\\
\label{eq:Rda_sopra}
\eps_0^{\frac{n}{2s}}\le \max\limits_j\int\limits_0^{2}\int\limits_{B_2'({x'_j})}
|u_h|^\sstar~dx'dx_1
\le \int\limits_0^{2}\int\limits_{B_2'(0)}
|u_h|^\sstar~dx'dx_1
\le  \left(2\eps_0\right)^{\frac{n}{2s}}.
\end{gather}
Up to a subsequence, we have that $u_h\to u$ weakly in  ${\widetilde{\mathcal D}^s(\R^n_+)}$.
To conclude the proof we show that $u\neq 0$.

Assume by contradiction that
$u_h\to 0$ weakly in  ${\mathcal D^s_{\rm R}(\R^n_+)}$. 
Ekeland's variational principle guarantees the existence of 
a sequence ${f}_h\to 0$ in the dual space  ${\mathcal D^s_{\rm R}(\R^n_+)}'$,
such that
\begin{equation}
\label{eq:Rtu_equation}
\DsNR u_h=
|u_h|^{\sstar-2}u_h+{f}_h\qquad\textrm{in ~$\mathcal D^s_{\rm R}(\R^n_+)'$}.
\end{equation}
Take $\f\in\mathcal C^\infty_0(-2,2)$ such that 
$\f\equiv 1$ on $(-1,1)$ and define $\f_j(x')=\f(|x'-x'_j|)$, $j=1,\dots,\tau$.
Then $\psi_j(x_1,x'):=\f(x_1)\f_j(x')$ has compact support in $(-2,2)\times B'_2(x'_j)$ and
$\psi_j\equiv 1$ on $B'_{1}(x_j)$ for any $j=1,\dots,\tau$. In addition,
$\psi_j^2u_h$ is a bounded sequence in $\mathcal D^s_{\rm R}(\R^n_+)$ by Lemma \ref{L:technicalR}.
We use $\psi_j^2u_h$ as test function
in (\ref{eq:Rtu_equation}) to find
\begin{equation}
\label{eq:Rmultline}
\langle\DsNR u_h,\psi_j^2u_h\rangle
=\int\limits_{\R^n_+} |u_h|^{\sstar-2}|\psi_j u_h|^2~dx+o(1)~\!.
\end{equation}
Thanks to H\"older inequality and (\ref{eq:Rda_sopra}) we can estimate
\begin{eqnarray}
\nonumber
\int\limits_{\R^n_+} \!|u_h|^{\sstar-2}|\psi_j u_h|^2dx
&\le& \Big(\int\limits_0^{2}\!\!\!\int\limits_{~\!B_2'(x'_j)}\!\!\!|u_h|^\sstar dx'dx_1\Big)^{\!\frac{2s}{n}}
\|\psi_ju_h\|^2_{L^\sstar(\R^n_+)}\\
&\le& 2\eps_0 \|\psi_ju_h\|^2_{L^\sstar(\R^n_+)}~\!.
\label{eq:numero}
\end{eqnarray}
We use Lemma \ref{L:technicalR}, Lemma \ref{L:spaceR} and the definition of $\mathcal S_R= \mathcal S_s^R(\R^n_+)$ to get
$$
\langle\DsNR u_h,\psi_j^2u_h\rangle = {E}(\psi_ju_h)+o(1)\ge {\mathcal S_R}\|\psi_ju_h\|^2_{L^\sstar(\R^n_+)}+o(1).
$$
Taking (\ref{eq:Rmultline}) into account, we arrive at
\begin{equation}
{\mathcal S_R}\|\psi_ju_h\|^2_{L^\sstar(\R^n_+)} \le 2\eps_0 \|\psi_ju_h\|^2_{L^\sstar(\R^n_+)}+o(1)~\!,
\label{eq:Rchain}
\end{equation}
that gives $\|\psi_ju_h\|_{L^\sstar(\R^n_+)}=o(1)$, because
$2\eps_0<{\mathcal S_R}$. 
Thus, using (\ref{eq:Rcover}) and recalling that $\psi_j\equiv 1$ on $(0,1)\times B'_1(x'_j)$, we obtain
$$
\int\limits_0^1\!\!\int\limits_{B'_2(0)}|u_h|^\sstar~\!dx'dx_1
\le \sum_{j=1}^\tau~\int\limits_0^1\!\!\int\limits_{B'_1(x_j)}|u_h|^\sstar~\!dx'dx_1\le
\sum_{j=1}^\tau~
\irnplus|\psi_ju_h|^\sstar dx %\int\limits_0^1\!\!\int\limits_{B'_1(x_j)}|\psi_ju_h|^\sstar~\!dx'dx_1
=o(1),
$$
that together with the first inequality in (\ref{eq:Rda_sopra}) gives
\begin{equation}
\label{eq:Rcappa}
\int\limits_1^2\!\!\int\limits_{B'_2(0)}|u_h|^\sstar~dx'dx_1\ge \eps_0^{\frac{n}{2s}}+o(1).
\end{equation}
Now we take a cut-off function $\phi\in\mathcal C^\infty_0(\R^n_+)$ such that
$\phi\equiv 1$ on $(1,2)\times B'_2(0)$. We test (\ref{eq:Rtu_equation}) with $\phi^2u_h\in\widetilde{\mathcal D}^s(\R^n_+)$ to get
\begin{equation}
\label{eq:spring}
\langle\DsNR u_h,\phi^2u_h\rangle
=\int\limits_{\R^n_+} |u_h|^{\sstar-2}|\phi u_h|^2~dx+o(1)~\!.
\end{equation}
Since $\text{supp}(\phi)\subset \R^n_+$, by 
\cite[Lemma 2.1]{MNcn} 
and thanks to the Sobolev inequality we obtain
$$
\langle\DsNR u_h,\phi^2u_h\rangle=\mathcal E_s(\phi u_h;\R^n\!\!\times\!\R^n)+o(1)
\ge \mathcal S_s\|\phi u_h\|_{L^\sstar(\R^n_+)}+o(1).
$$
Therefore, estimating the right hand side of (\ref{eq:spring}) via H\"older inequality we obtain
\begin{equation}
\label{eq:miami}
\mathcal S_s\|\phi u_h\|_{L^\sstar(\R^n_+)}\le 
\|u_h\|^{\sstar-2}_{L^\sstar(\R^n_+)}\|\phi u_h\|^2_{L^\sstar(\R^n_+)}+o(1)
={\mathcal S_R}\|\phi u_h\|^2_{L^\sstar(\R^n_+)}+o(1).
\end{equation}
Now we recall that $\mathcal S_R<\mathcal S_s$ and
$\phi\equiv 1$ on $(1,2)\times B'_2(0)$. Thus (\ref{eq:miami}) gives  
$$
\int\limits_1^2\!\!\int\limits_{B'_2(0)}|u_h|^\sstar~\!dx_1dx'=o(1)~\!.$$ 
We reached a contradiction with (\ref{eq:Rcappa}),
that concludes the proof.
\QED

\begin{Remark}[Euler-Lagrange equations]
\label{R:REL}
Any extremal for 
$S^{\rm R}_s(\R^n_+)$ solves, up to a Lagrange multiplier, the nonlocal differential equation
$\DsNR u=|u|^{\sstar-2}u$ in $\R^n_+$.
Standard arguments and \cite[Remark 2.5]{MNmp} imply that $u$ has constant sign on $\R^n_+$.
We can assume that $u$ is nonnegative on $\R^n$, so that $u$ is a weak solution to
$$
\DsNR u=u^{\sstar-1}\qquad\text{in $\R^n_+$.}
$$
Thus $u$ is lower semicontinuous and positive by 
the strong maximum principle in  \cite[Corollary 4.3]{MNmp}. 

Next we deal with boundary conditions.
First assume 
$0<s<\frac12$. Arguing as in \cite[Sec. 2.10.2]{Tr} one gets that
${\mathcal D}^s_{\rm R}(\R^n_+)=\widetilde{\mathcal D}^s(\R^n_+)$ with equivalent norms. Further, formula  (\ref{eq:N_sum}) gives
\begin{equation}
\label{eq:problem1,5}
S^{\rm R}_s(\R^n_+)=\inf_{v\in \widetilde{\mathcal D}^s(\R^n_+)\atop\scriptstyle  v\ne 0}
\frac{\langle\Ds v,v\rangle-\gamma_s\|x_1^{-s}v\|^2_{L^2(\R^n_+)}}{\|v\|^2_{L^\sstar(\R^n_+)}}~\!,
\end{equation}
where $\gamma_s$ has been defined in (\ref{eq:gamma}), and the
 minimization problems
(\ref{eq:problem1,5}), (\ref{eq:problemR}) are  equivalent. Thus 
$u$ solves the Dirichlet's problem
$$
\Ds u=\gamma_s|x_1|^{-2s}u+u^{\sstar-1}\qquad\text{in $\R^n_+$,~~$u\equiv 0$~in $\R^n_-$.}
$$
We cite the papers \cite{DLV, Frank, MNcn} for related results.
 
If $s=\frac 12$ it is not clear whether one can even talk about
boundary conditions for $u$ (however, following the arguments in \cite[Sec. 4.3.2]{Tr} one can see that
$u$ can be approximated in $\mathcal D^\frac12_{\rm R}(\R^n_+)$ by
a sequence of functions in $C^\infty_0(\R^n_+)$).

For $s\in(\frac12,1)$ we can use the results in \cite{G1, G2, W2} to conclude that 
$u$ satisfies the Neumann-type boundary condition
$$
\mathcal N_s u(x'):=-(2s-1)\lim_{x_1\to 0^+}x_1^{1-2s}(u(x_1,x')-u(0,x'))=0,\qquad x'\in\partial \R^n_+.
$$
\end{Remark}

\section{Neumann Spectral Laplacian: proof of Theorem \ref{T:Spectral}}
\label{S:Spectral}

The {\em Neumann Spectral  fractional Laplacian} is the $s$-th power of standard Neumann Laplacian in the sense of spectral theory.
For the Laplacian in $\R^n_+$ this gives representation
$$
\langle\DsNSp u,u\rangle=\irnplus |\xi|^{2s}|{\cal F}u(\xi)|^2\,d\xi,
$$
where
$$
{\cal F}u(\xi_1,\xi')=\frac 2{(2\pi)^{\frac n2}}\irnplus \cos(x_1\xi_1)\exp(-ix'\cdot\xi')u(x)\,dx
$$
is the cosine Fourier transform.

Denote by $\hat u$ the even extension of $u$, see (\ref{eq:even_ex}). It is easy to see that
$$
\langle\DsNSp u,u\rangle=\frac 12 \langle\Ds \hat u,\hat u\rangle~\!,\qquad \|u\|^2_{L^\sstar(\R^n_+)}
=\frac 1{2^{1-\frac {2s}n}}\|\hat u\|^2_{L^\sstar(\R^n)}.
$$
Thus $\langle\DsNSp u,u\rangle$ is finite just for $u$ in the space ${\mathcal D}^s_{\rm R}(\R^n_+)$.
Moreover
$$
\mathcal S^{\rm Sp}_s(\R^n_+):=\inf_{u\in {\mathcal D}^s_{\rm R}(\R^n_+)\atop\scriptstyle  u\ne 0}
\frac{\langle\DsNSp u,u\rangle}{\|u\|^2_{L^\sstar(\R^n_+)}}\ge 2^{-\frac {2s}n}\mathcal S_s~\!,
$$
so the function $U_{\!s}$ defined in (\ref{U}) easily provides the value $2^{-\frac {2s}n}\mathcal S_s$,
and Theorem \ref{T:Spectral} is proved.
\QED

\section{Semirestricted Laplacian and proof of Theorem \ref{T:Semirestricted}} 
\label{S:Semirestricted}

We start with few  remarks about the space 
${\mathcal D}^s_{\rm Sr}(\R^n_+)$. 

\begin{enumerate}
\item
First of all we have 
\begin{equation}
\label{eq:all}
\mathcal E_s(u;\R^{2n}\setminus(\R^n_-)^2)=
\mathcal E_s(u;\R^n_+\!\times\!\R^n_+)+2\mathcal E_s(u;\R^n_+\!\times\!\R^n_-)
\quad\text{for any $u\in {\mathcal D}^s_{\rm Sr}(\R^n_+)$.}
\end{equation}

\item
\label{I:D}
 $\mathcal D^s(\R^n)\subset \mathcal D^s_{\rm Sr}(\R^n)$ and 
\begin{equation}
\label{eq:few}
\mathcal E_s(u;\R^{2n}\setminus(\R^n_-)^2)=\mathcal E_s(u;\R^n\!\times\R^n)-\mathcal E_s(u;\R^n_-\!\times\R^n_-)
\quad\text{for any $u\in {\mathcal D}^s(\R^n)$.} 
\end{equation}
In particular $\widetilde{\mathcal D}^s(\R_+^n)=\{u\in {\mathcal D}^s_{\rm Sr}(\R^n_+)~|~ u\equiv 0~~
\text{on $\R^n_+$}\}$, and 
$$
\mathcal E_s(u;\R^{2n}\setminus(\R^n_-)^2)=\mathcal E_s(u;\R^n\!\times\R^n)\quad \text{
for any $u\in \widetilde{\mathcal D}^s(\R_+^n)$.}
$$

\item
\label{I:3}
Let
$u\in {\mathcal D}^s_{\rm Sr}(\R^n_+)$ and assume  $u=0$ a.e. on $\R^n_+$.  From 
(\ref{eq:all}) it follows that $u\in L^2(\R^n;|x_1|^{-2s}dx)$, that is the space of  functions on 
$\R^n$ that  are square integrable with respect 
to the measure $|x_1|^{-2s}dx$. More precisely
\begin{equation}
\label{eq:vanishing}
\mathcal E_s(u;\R^{2n}\setminus(\R^n_-)^2)=\gamma_s\displaystyle\int\limits_{\R^n_-}|x_1|^{-2s}|u|^2~\!dx~\!,
\end{equation}
where $\gamma_s$ is defined in (\ref{eq:gamma}). Also the converse is true, namely, if $u\in L^2(\R^n;|x_1|^{-2s}dx)$ and $u=0$ on $\R^n_+$, then
$u\in {\mathcal D}^s_{\rm Sr}(\R^n_+)$ and (\ref{eq:vanishing}) holds.
In particular, ${\mathcal D}^s(\R^n)$ is properly contained in ${\mathcal D}^s_{\rm Sr}(\R^n_+)$.

\item
Let $u\in {\mathcal D}^s_{\rm Sr}(\R^n_+)$. Clearly $u|_{\R^n_+}\in \mathcal D^s_{\rm R}(\R^n_+)$. 
We decompose $u$ via the even extension operator in (\ref{eq:even_ex}). Precisely we write
$$
u=\widehat{u|_{\R^n_+}}+\big(u-\widehat{u|_{\R^n_+}}\big).
$$
We have $\widehat{u|_{\R^n_+}}\in {\mathcal D}^s(\R^n)\subset {\mathcal D}^s_{\rm Sr}(\R^n_+)$,
hence $u-\widehat{u|_{\R^n_+}} \in {\mathcal D}^s_{\rm Sr}(\R^n_+)$. Thus
$u-\widehat{u|_{\R^n_+}}\in L^2(\R^n;|x_1|^{-2s}dx)$ by Item \ref{I:3}. We have shown 
$\mathcal D^s_{\rm Sr}(\R^n_+)\subset\mathcal D^s(\R^n)+L^2(\R^n;|x_1|^{-2s}dx)$.
Actually, it is easy to check that 
$$
{\mathcal D}^s_{\rm Sr}(\R^n_+)=\big\{v+\chi_{\R^n_-}g~|~ v\in \mathcal D^s(\R^n)~,~~g\in L^2(\R^n;|x_1|^{-2s}dx)~\big\}.
$$
\end{enumerate}

\begin{Lemma}
\label{L:Sr_space}
The space ${\mathcal D}^s_{\rm Sr}(\R^n_+)$ inherits an Hilbertian structure from the norm
$$
\|u\|_{ {\mathcal D}^s_{\rm Sr}(\R^n_+)}^2=\mathcal E_s(u;\R^{2n}\setminus(\R^n_-)^2).
$$
Moreover, the restriction operator $u\mapsto u|_{\R^n_+}$ is continuous ${\mathcal D}^s_{\rm Sr}(\R^n_+)\to 
{\mathcal D}^s_{\rm R}(\R^n_+)$ and ${\mathcal D}^s_{\rm Sr}(\R^n_+)\to L^\sstar(\R^n_+)$.
\end{Lemma}

\proof
Let $u\in {\mathcal D}^s_{\rm Sr}(\R^n_+)$.  If 
$\mathcal E_s(u;\R^{2n}\setminus(\R^n_-)^2)=0$ then $u\equiv 0$ on $\R^n$ by (\ref{eq:all}) and (\ref{eq:vanishing}).
Thus $\|u\|_{ {\mathcal D}^s_{\rm Sr}(\R^n_+)}$ is a norm on ${\mathcal D}^s_{\rm Sr}(\R^n_+)$.

The conclusions about the restriction operator $u\mapsto u|_{\R^n_+}$ are immediate, use
also the continuity of the embedding ${\mathcal D}^s_{\rm R}(\R^n_+)\hookrightarrow L^\sstar(\R^n_+)$
given by Lemma \ref{L:spaceR}.

To check completeness one can adapt the argument for \cite[Proposition 3.1]{DRoV} or argue as follows. 
Let $u_h$ be a Cauchy sequence in ${\mathcal D}^s_{\rm Sr}(\R^n_+)$. Then 
${u_h}|_{\R^n_+}$ is a Cauchy sequence in ${\mathcal D}^s_{\rm R}(\R^n_+)$. We write $u_h$ as
$$
u_h=w_h+(u_h-w_h)~,\quad w_h= \widehat{u|_{\R^n_+}}\in \mathcal D^s(\R^n)~,\quad
u_h-w_h\in L^2(\R^n;|x_1|^{2s}dx),
$$
and $u_h-w_h\equiv 0$ on $\R^n_+$. Then $w_h$ is a Cauchy sequence in $\mathcal D^s(\R^n)$
and in $\mathcal D^s_{\rm Sr}(\R^n_+)$, see (\ref{eq:numero0}), and
$u_h-w_h$ is a Cauchy sequence in $L^2(\R^n;|x_1|^{2s}dx)$ by (\ref{eq:vanishing}). Thus $w_h\to w$  in 
$\mathcal D^s(\R^n)$, $u_h-w_h\to g$ in $L^2(\R^n;|x_1|^{2s}dx)$ and therefore
$u_h\to w+g$ in ${\mathcal D}^s_{\rm Sr}(\R^n_+)$.
\QED

For $u\in \mathcal D^s_{\rm Sr}(\R^n_+)$ we introduce 
the distribution 
$\DsNR u\in \mathcal D^s_{\rm Sr}(\R^n_+)'$ by
\begin{equation}
\label{eq:Sr_diff}
\langle \DsNS u,v\rangle= \frac{C_{n,s}}{2}\iint\limits_{\R^{2n}\setminus(\R^n_-)^2} 
\frac{(u(x)-u(y))(v(x)-v(y))}{|x-y|^{n+2s}}~\!\!dxdy~,\quad
v\in \mathcal D^s_{\rm Sr}(\R^n_+).
\end{equation}

\medskip
Before going further, let us try to explain why we need more preliminary results to prove 
Theorem \ref{T:Semirestricted}.

Any bounded minimizing sequence $u_h$ for (\ref{eq:DRoV})
has a subsequence $u_h$ such that $u_h\to u$ weakly in ${\mathcal D}^s_{\rm Sr}(\R^n_+)$. 
We surely have
$\mathcal E_s(u;\R^{2n}\setminus(\R^n_-)^2)\le {\mathcal S}^{\rm Sr}_s(\R^n_+)$; moreover 
we can say that ${u_h}|_{\R^n_+}\to u|_{\R^n_+}$ weakly in $\mathcal D^s_{\rm Sr}(\R^n_+)$,
weakly in $L^\sstar(\R^n_+)$, strongly in $L^p_{\rm loc}(\overline{\R}^n_+)$ for $p\in[1,\sstar)$ and
pointwise almost everywhere on $\R^n_+$. 
However, no informations on the pointwise (for instance) convergence of $u_h$ on $\R^n_-$ are available, and
we can not go further with the study of the behavior of $u_h$ on $\R^n$.
In essence, to overcome this technical difficulty we  
move from (\ref{eq:DRoV}) to an equivalent minimization problem
that inherits better (local) compactness properties on $\R^n$, being  settled on a
smaller function space. 

The first step consists in finding the ''best extension'' $P_{\!s}u \in \mathcal D^s_{\rm Sr}(\R^n_+)$
of $u|_{\R^n_+}$, 
for any function $u\in \mathcal D^s_{\rm Sr}(\R^n_+)$  (see also \cite[Section 5]{DRoV}). 
Recall that the value of the constant $\gamma_s$ is given in (\ref{eq:gamma}).
The three  lemmata that follow are proved in Appendix \ref{S:P}. 

\begin{Lemma}
\label{L:extension}
\begin{itemize}
\item[$i)$]
Let $u\in \mathcal D^s_{\rm Sr}(\R^n_+)$. The function
\begin{equation}
\label{eq:P_explicit}
(P_{\!s}u)(x)=\begin{cases}u(x)&\text{if $x\in\R^n_+$}\\
\displaystyle{{\frac{C_{n,s}}{\gamma_s}}
|x_1|^{2s}\displaystyle\irnplus\frac{u(y)}{|x-y|^{n+2s}}~\!dy}&\text{if $x\in\R^n_-$.}
\end{cases}
\end{equation}
is the unique solution to the convex
minimization problem
\begin{equation}
\label{eq:P_minimization}
\min_{\scriptstyle\omega\in \mathcal D^s_{\rm Sr}(\R^n_+)\atop 
\scriptstyle\omega|_{\R^n_+}=u} \mathcal E_s(\omega;\R^{2n}\setminus(\R^n_-)^2);
\end{equation}
\item[$ii)$] The linear operator $P_{\!s}$ is orthoprojector in $\mathcal D^s_{\rm Sr}(\R^n_+)$
that is, $P_{\!s}^2=P_{\!s}$ and $P_{\!s}^*=P_{\!s}$;
\item[$iii)$] If $u\in \mathcal D^s_{\rm Sr}(\R^n_+)$ then
\begin{equation}
\label{eq:orthogonal}
\mathcal E_s(u-P_{\!s}u;\R^{2n}\setminus(\R^n_-)^2)=\mathcal E_s(u;\R^{2n}\setminus(\R^n_-)^2)-\mathcal E_s(P_{\!s}u;\R^{2n}\setminus(\R^n_-)^2)~\!.
\end{equation}
\end{itemize}
\end{Lemma}
Now we study the image of the operator $P_{\!s}$. We put 
$$
\mathcal R^s_{\rm Sr}(\R^n_+):=Im(P_{\!s})
= \{u\in \mathcal D^s_{\rm Sr}(\R^n_+)~|~u=P_{\!s}u\}~\!.
$$
Obviously, $\mathcal R^s_{\rm Sr}(\R^n_+)$ is a Hilbert subspace of $\mathcal D^s_{\rm Sr}(\R^n_+)$.

\begin{Lemma}
\label{L:Pcompactness}
\begin{itemize}
\item[$i)$] $\mathcal R^s_{\rm Sr}(\R^n_+)$ is continuously embedded into $L^\sstar(\R^n)$;
\item[$ii)$] $\mathcal R^s_{\rm Sr}(\R^n_+)$ is compactly embedded into $L^p_{\rm loc}(\R^n)$,
for any $p\in[1,\sstar)$. That is, for any sequence $u_h\in \mathcal R^s_{\rm Sr}(\R^n_+)$ such that $u_h\to 0$ weakly in 
$\mathcal D^s_{\rm Sr}(\R^n_+)$,
there exists a subsequence $u_h$ such that $u_h\to 0$ in $L^p_{\rm loc}(\R^n)$.
\end{itemize}
\end{Lemma}

\begin{Lemma}
\label{L:differential}
Let $u\in \mathcal R^s_{\rm Sr}(\R^n_+)$. Then
\begin{itemize}
\item[$i)$] $\Ds u$ is a distribution on $\R^n$, and $\DsNS u= \Ds u$ in $\mathcal D'(\R^n_+)$;
\item[$ii)$] $\DsNS u=\DsNR u=0$ in $\mathcal D'(\R^n_-)$ and almost everywhere on $\R^n_-$. 
\end{itemize}
\end{Lemma}

\begin{Remark}
\label{R:equivalent}
The Sobolev constant ${\mathcal S}^{\rm Sr}_s(\R^n_+)$ coincides with
\begin{equation}
\label{eq:min_equivalent}
\tilde {\mathcal S}^{\rm Sr}_s(\R^n_+)=\inf_{u \in {\mathcal R}^s_{\rm Sr}(\R^n)\atop\scriptstyle  u\ne 0}
\frac{\mathcal E_s(u;\R^{2n}\setminus(\R^n_-)^2)}{\|u\|^2_{L^\sstar(\R^n_+)}}.
\end{equation}
Moreover, the minimization problems in (\ref{eq:DRoV}), (\ref{eq:min_equivalent}) are equivalent, 
that is, $u$ achieves ${\mathcal S}^{\rm Sr}_s(\R^n_+)$ if and only if $u=P_{\!s}u$ and $u$ achieves
$\tilde {\mathcal S}^{\rm Sr}_s(\R^n_+)$.
\end{Remark}

The following statement is the analog of Lemma \ref{L:technicalR} for the Semirestricted Laplacian.

\begin{Lemma}
\label{L:technicalSr}
Let $u\in \mathcal R^s_{\rm Sr}(\R^n_+)$,   $\f\in\mathcal C^\infty_0(\R^n)$, and let $\Omega$ be a bounded
domain containing the support of $\f$. Then $\f u\in \mathcal D^s_{\rm Sr}(\R^n_+)$
and
\begin{gather*}
\mathcal E_s(\f u;\R^{2n}\setminus(\R^n_-)^2)\le c(\f)\big(\mathcal E_s(u;\R^{2n}\setminus(\R^n_-)^2)+
\|u\|_{L^2(\Omega)}^2\big);\\
\big|\mathcal E_s(\f u;\R^{2n}\setminus(\R^n_-)^2)\!-\!\langle \DsNS u,\f^2 u\rangle\big|\!\!\le\! c(\f) 
\|u\|_{L^2(\Omega)}\!\!\left(
\|u\|_{L^2(\Omega)}
+\mathcal E_s(u;\R^{2n}\setminus(\R^n_-)^2)^\frac12\!\right)\!.
\end{gather*}
\end{Lemma}

\proof 
We keep the notation in (\ref{eq:to_simplify}). First, recall that $u\in L^2_{\rm loc}(\R^n)$ by Lemma \ref{L:Pcompactness}. 
Then we estimate
\begin{eqnarray*}
\mathcal E_s(\f u;\Omega^+\!\!\times\!\R^n_-)&=&c\iint\limits_{\Omega^+\times\R^n_-}
\frac{\big((\f u)(x)-(\f u)(y)\pm\f(x)u(y)\big)^2}{|x-y|^{n+2s}}~\!dxdy\\
&\le& c(\|\f\|_\infty)\mathcal E_s(u;\R^{2n}\setminus(\R^n_-)^2)
+\int\limits_{\Omega^+} |u(x)|^2\big(\irn\Psi_\f(x,y)~\!dy\big)dx~\!.
\end{eqnarray*}
Thus (\ref{eq:Psi_estimate}) gives
$\mathcal E_s(\f u;\Omega^+\!\!\times\!\R^n_-)\le c(\f)
\big(\mathcal E_s(u;\R^{2n}\setminus(\R^n_-)^2)+ \|u\|_{L^2(\Omega)}^2\big)$. Arguing
in a similar way one finds the same bound for 
$\mathcal E_s(\f u;\R^n_+\times\Omega^-)$. 
Since $(\f u)(x)-(\f u)(y)= 0$ unless $x\in\Omega$ or $y\in\Omega$, we infer that
$$
\mathcal E_s(\f u;\R^n_+\!\!\times\!\R^n_-)\le 
c(\f)\big(\mathcal E_s(u;\R^{2n}\setminus(\R^n_-)^2)+ \|u\|_{L^2(\Omega)}^2\big)~\!.
$$
The first inequality in Lemma \ref{L:technicalSr} follows by  Lemma \ref{L:technicalR} and by the equality
$$
\mathcal E_s(\f u;\R^{2n}\setminus(\R^n_-)^2)=
\mathcal E_s(\f u;\R^n_+\!\times\!\R^n_+)+2\mathcal E_s(\f u;\R^n_+\!\times\!\R^n_-)~\!.
$$
Next, we compute
$$
\mathcal E_s(\f u;\R^{2n}\setminus(\R^n_-)^2)- \langle \DsNS u,\f^2 u\rangle=
\frac{C_{n,s}}{2}\iint\limits_{\R^{2n}\setminus(\R^n_-)^2}{u}(x){u}(y)\Psi_\f(x,y)~\!dxdy=:\mathcal B_\f.
$$
Since $\Psi_\f(x,y)\equiv 0$ on $(\R^n\setminus\Omega)^2$ and since
$\big(\R^{2n}\setminus(\R^n_-)^2\big)\setminus(\R^n\setminus\Omega)^2$ is contained in the union of
$\Omega\times\Omega$ with the sets
$$
A:=\big[ \Omega\times(\R^n_+\setminus\Omega)\big]\cup\big[\Omega^+\times(\R^n_-\setminus\Omega)\big]~,~~
\hat A=\big[(\R^n_+\setminus\Omega)\times\Omega\big]\cup\big[(\R^n_-\setminus\Omega)\times \Omega^+\big],
$$
we can estimate
$$
\displaystyle{
c~\!|\mathcal B_\f|\le \iint\limits_{\Omega\!\times\Omega} |u(x)u(y)|\Psi_\f~\!dxdy+
\iint\limits_{A}|u(x)|\frac{|u(y)||\f(x)|^2}{|x-y|^{n+2s}}dxdy}.
$$  
We put 
\begin{gather*}
J_1=\iint\limits_{\Omega\times\Omega} |u(x)u(y)|\Psi_\f~\!dxdy~,\quad
J_2=\iint\limits_A \frac{|u(x)\f(x)|^2}{|x-y|^{n+2s}}~\!dxdy\\
J_3=\iint\limits_{A} \frac{|u(x)-u(y)|~\!|u(x)\f(x)|}{|x-y|^{n+2s}}~\!dxdy.
\end{gather*}
Arguing as in the proof of Lemma \ref{L:technicalR} one gets
$c|\mathcal B_\f|  \le J_1 +J_2+\|\f\|_\infty J_3$ and estimates
$$
J_1, J_2\le c(\f)\int\limits_{\Omega}|u(x)|^2~\!dx~,~~
J_3^2\le 2\mathcal E_s(u;A) J_2 %\int\limits_{\Omega}|u(x)|^2dx
\le c(\f)\mathcal E_s(u;\R^{2n}\setminus(\R^n_-)^2) \int\limits_{\Omega}|u(x)|^2dx,
$$
that concludes the proof.
\QED

\bigskip
\noindent
{\bf Proof of Theorem \ref{T:Semirestricted}.}
The first claim follows from
$$
\mathcal S^{\rm Sr}_s(\R^n_+)\le \inf_{u\in\mathcal C^\infty_0(\R^n_+)\atop u\neq 0}
\frac{\mathcal E_s(u;\R^{2n}\setminus(\R^n_-)^2)}{\|u\|^2_{L^\sstar(\R^n_+)}}=
\inf_{u\in\mathcal C^\infty_0(\R^n_+)\atop u\neq 0}
\frac{\mathcal E_s(u;\R^{n}\times\R^n)}{\|u\|^2_{L^\sstar(\R^n)}}= \mathcal S_s.
$$
Now assume ${\mathcal S}^{\rm Sr}_s(\R^n_+)<\mathcal S_s$. We have to show that
there exists a minimizer for ${\mathcal S}^{\rm Sr}_s(\R^n_+)$.
The argument does not differ too much from that used in the proof of 
Theorem \ref{T:Restricted}; we limit ourselves to point out the main changes. 
With respect to the Restricted case, the main differences concern the role played 
by the operator $P_{\!s}:\mathcal D^s_{\rm Sr}(\R^n_+)\to \mathcal D^s_{\rm Sr}(\R^n_+)$.

As in the proof of Theorem \ref{T:Restricted} we restrict ourselves to the case $n\ge 2$ and we put 
$$\mathcal S_{\rm Sr}={\mathcal S}^{\rm Sr}_s(\R^n_+)~,\quad E(u)=\mathcal E_s(u;\R^{2n}\setminus(\R^n_-)^2).
$$
We only need to exhibit a minimizing sequence for ${\mathcal S}^{\rm Sr}_s(\R^n_+)$
that weakly converges to a nontrivial limit. Fix a number $\eps_0\in (0,\frac12 \mathcal S_{\rm Sr})$. 
Argue as in the proof of Theorem \ref{T:Restricted} to construct 
a minimizing sequence $u_h\in\mathcal D^s_{\rm Sr}(\R^n_+)$ satisfying
\begin{equation}
\label{eq:gemma}
\|{u}_h\|^\sstar_{L^\sstar(\R^n_+)}=\mathcal S_{\rm Sr}^{\frac{n}{2s}}~,~~
{E}({u}_h)=\mathcal S_{\rm Sr}^{\frac{n}{2s}}+o(1)
\end{equation}
and such that (\ref{eq:Rda_sopra}) holds. %for any ${z}\in\R^{n-1}$.
We can assume that ${u}_h\to {u}$ weakly in  ${\mathcal D^s_{\rm Sr}(\R^n_+)}$. 
By contradiction suppose that $u=0$.
Consider now the sequence $P_{\!s}u_h\in \mathcal R^s_{\rm Sr}(\R^n_+)$, that is bounded in
$\mathcal R^s_{\rm Sr}(\R^n_+)$ and  satisfies
$$
\|{P_{\!s}u}_h\|^\sstar_{L^\sstar(\R^n_+)}=\|{u}_h\|^\sstar_{L^\sstar(\R^n_+)}=\mathcal S_{\rm Sr}^{\frac{n}{2s}}~,~~
\mathcal S_{\rm Sr}^{\frac{n}{2s}}\le {E}({P_{\!s}u}_h)\le {E}({u}_h)=\mathcal S_{\rm Sr}^{\frac{n}{2s}}+o(1).
$$
In particular, $P_{\!s}u_h$ is a minimizing sequence for (\ref{eq:DRoV}) (and for the equivalent minimization problem
(\ref{eq:min_equivalent})). 

Next we notice that $P_{\!s}u_h-u_h\to 0$ in $\mathcal D^s_{\rm Sr}(\R^n_+)$
by (\ref{eq:orthogonal}). 
We infer that
$P_{\!s}u_h\to 0$
weakly in $\mathcal  R^s_{\rm Sr}(\R^n_+)$. 
In other words, $P_{\!s}u_h$ enjoys the same properties as the sequence $u_h$
(including (\ref{eq:Rda_sopra}), as $P_{\!s}u_h\equiv u_h$ on $\R^n_+$), and in addition
$P_{\!s}u_h\in \mathcal R^s_{\rm Sr}(\R^n_+)$. In order to simplify notation, from now on
we write  $u_h=P_{\!s}u_h$.

By Ekeland's variational principle we can assume that
there exists a sequence ${f}_h\to 0$ in  ${\mathcal D^s_{\rm Sr}(\R^n_+)}'$,
such that
\begin{equation}
\label{eq:Stu_equation}
\DsNS {u}_h=
\chi_{\R^n_+}|{u}_h|^{\sstar-2}{u}_h+f_h\qquad\textrm{in $\mathcal D^s_{\rm Sr}(\R^n_+)'$}.\\
\end{equation}
Take points 
$x'_1,\cdots x'_\tau\in \R^{n-1}$ and 
cut-off functions $\psi_j$, $j=1,\dots,\tau$
as in the proof of Theorem \ref{T:Restricted} and test (\ref{eq:Stu_equation}) with
$\psi_j^2u_h\in \mathcal D^s_{\rm Sr}(\R^n_+)$.
Use Lemma \ref{L:technicalSr}, the last inequality in (\ref{eq:Rda_sopra}) and
adapt the computations for (\ref{eq:Rchain}) to obtain
$$
\mathcal S_{\rm Sr}\|\psi_j{u}_h\|^2_{L^\sstar(\R^n_+)}\le
2\eps_0 \|\psi_j{u}_h\|^2_{L^\sstar(\R^n_+)}+o(1)~\!.
$$
Thus 
(\ref{eq:Rcappa}) holds also in this case, because $2\eps_0<\mathcal S_{\rm Sr}$. Next we
take $\phi\in\mathcal C^\infty_0(\R^n_+)$ such that
$\phi\equiv 1$ on $(1,2)\times B'_2(0)$. 
Notice that $\phi u_h\in \widetilde{\mathcal D}^s(\R^n_+)\subset {\mathcal D}^s(\R^n)$
and $\phi u_h\to 0$ in $L^2(\R^n)$. In particular, from
Lemma \ref{L:technicalSr} and thanks to the Sobolev inequality
we get
$$
\langle \DsNS {u}_h,\phi^2u_h\rangle=
E(\phi u_h)+o(1)=\mathcal E_s(\phi u_h;\R^n\!\times\!\R^n)+o(1)
\ge \mathcal S_s \|\phi u\|^2_{L^\sstar(\R^n_+)}+o(1).
$$
Therefore,  testing  (\ref{eq:Stu_equation}) with 
$\phi^2u_h$
and  using H\"older inequality  we obtain
$$
\mathcal S_s\|\phi{u}_h\|^2_{L^\sstar(\R^n_+)}\le\|{u}_h\|^{\sstar-2}_{L^\sstar(\R^n_+)}\|\phi{u}_h\|^2_{L^\sstar(\R^n_+)}+o(1)
=\mathcal S_{\rm Sr}\|\phi{u}_h\|^2_{L^\sstar(\R^n_+)}+o(1).
$$
We infer that $\|\phi{u}_h\|^2_{L^\sstar(\R^n_+)}=o(1)$, because $\mathcal S_{\rm Sr}<\mathcal S_s$. We reached a contradiction 
with (\ref{eq:Rcappa}), as $\phi\equiv 1$ on $(1,2)\times B'_2(0)$. 
\QED

\begin{Theorem}
\label{T:strictSr}
It holds that $\mathcal S^{\rm Sr}_s(\R^n_+)<\mathcal S_s$ provided that one of the following conditions is satisfied:
\begin{itemize}
\item[$i)$] $n\ge 2$ and $s$ is close enough to $1^-$;
\item[$ii)$] $n=1$ and $s$ is close enough to $\frac12^-$. 
\end{itemize}
\end{Theorem}

\proof
First, assume  $n\ge 3$. 
Recall that the function $U(x)=(1+|x|^2)^{\frac{2-n}{2}}\in\mathcal D^1(\R^n)$ achieves
the Sobolev constant in (\ref{eq:Sob1}).
For any $\eps>0$, we can find
a radial function $u_\eps\in\mathcal C^\infty_0(\R^n)$ that is close to $U$ in the $\mathcal D^1(\R^n)$ topology,
and such that
$$
\frac{\|\nabla u_\eps\|^2_{L^2(\R^n)}}{\|u_\eps\|^2_{L^{2^*}(\R^n)}}<\mathcal S+\eps.
$$
Our aim is to estimate from above $\mathcal S^{\rm Sr}_s(\R^n_+)$ via ${u_\eps}$, for $s\to 1^-$.
Clearly $|u_\eps|^\sstar\to |u_\eps|^{2^*}$ in $L^1(\R^n)$. It holds that
$\mathcal E_s(u_\eps;\R^{2n}\setminus(\R^n_-)^2)\to \frac12\|\nabla u_\eps\|^2_{L^2(\R^n)}$
as $s\to 1^-$, see Theorem \ref{P:s_to1} in Appendix \ref{A:A}.
Thus
$$
\limsup_{s\to 1^-} \mathcal S^{\rm Sr}_s(\R^n_+)\le 
\lim_{s\to 1^-}\frac{\mathcal E_s({u_\eps};\R^{2n}\setminus(\R^n_-)^2)}{\|{u_\eps}\|^2_{L^\sstar(\R^n_+)}}
= \frac{1}{2^{\frac2n}}\frac{\|\nabla u_\eps\|^2_{L^2(\R^n)}}{\|u_\eps\|^2_{L^{2^*}(\R^n)}}
<\frac{1}{2^{\frac2n}}(\mathcal S+\eps).
$$
Since $\mathcal S_s$ depends continuously on $s\in(0,n/2)$, see (\ref{eq:CoTa}),
for $\eps$ small enough and $s$ close enough to $1$
we see that $\mathcal S^{\rm Sr}_s(\R^n_+)<\mathcal S_s$.

Next, assume $n=1$ or $n=2$. We test $\mathcal S^{\rm Sr}_s(\R^n_+)$ with the function $U_{\!s}$ in (\ref{U}). Since
\begin{gather*}
\mathcal E_s({U_{\!s}};\R^n\!\!\times\!\R^n)=2\mathcal E_s({U_{\!s}};\R^n_+\!\!\times\!\R^n_+)+2\mathcal E_s({U_{\!s}};\R^n_+\!\!\times\!\R^n_-)
<4\mathcal E_s({U_{\!s}};\R^n_+\!\!\times\!\R^n_+)\\
\mathcal E_s({U_{\!s}};\R^{2n}\setminus(\R^n_-)^2)=\mathcal E_s({U_{\!s}};\R^n\!\!\times\!\R^n)-\mathcal E_s({U_{\!s}};\R^n_+\!\!\times\!\R^n_+)
\end{gather*}
we have $\mathcal E_s({U_{\!s}};\R^n\!\!\times\!\R^n)>\frac43\mathcal E_s({U_{\!s}};\R^{2n}\setminus(\R^n_-)^2) $, hence
$$
\mathcal S_s=\frac{\mathcal E_s({U_{\!s}};\R^n\!\times\!\R^n)}{\|{U_{\!s}}\|^2_{L^\sstar(\R^n)}}>
2^{\frac{2s}{n}}\frac{2}{3}~\!\frac{\mathcal E_s({U_{\!s}};\R^{2n}\setminus(\R^n_-)^2)}{\|{U_{\!s}}\|^2_{L^\sstar(\R^n_+)}}
\ge 2^{\frac{2s}{n}}\frac{2}{3}S^{\rm Sr}_s(\R^n_+).
$$
In particular, if $1>\frac{2s}{n}\ge \frac{\ln(3/2)}{\ln 2}$ then the desired strict inequality holds. This
condition is satisfied provided that $s$ is close enough to $1$ if $n=2$, and close enough to $\frac12$ if $n=1$.
\QED

\bigskip

\noindent
{\bf Conjecture}
{\em Quite likely, it holds that $\mathcal S^{\rm Sr}_s(\R^n_+)<\mathcal S_s$ for all admissible exponents $s$.}

\begin{Remark}[Euler-Lagrange equations]
Let $u\in {\mathcal D}^s_{\rm Sr}(\R^n_+)$ be an extremal for  $S^{\rm R}_s(\R^n_+)$. 
Then
 $u$ solves $\DsNS u=\chi_{\R^n_+}|u|^{\sstar-2}u$ in $\R^n$, up to a Lagrange multiplier.
 
 By Remark \ref{R:equivalent} and Lemma \ref{L:differential} 
 we have that $u\in {\mathcal R}^s_{\rm Sr}(\R^n_+)$ and $u$ solves
$$
\Ds u=|u|^{\sstar-2}u\quad\text{ in $\R^n_+$},\quad \DsNR u=0\quad\text{ in $\R^n_-$.}
$$
Using standard arguments and \cite[Remark 2.5]{MNmp}, one can easily prove 
that $u$ can not change sign on $\R^n_+$, so that we can assume that $u$ is nonnegative on $\R^n_+$.
By \cite[Corollary 4.4]{MNmp}  a strong maximum principle holds for the operator $\DsNS$.
In particular, we  have that $u$ is lower semicontinuous and positive in $\R^n_+$.
\end{Remark}

\section{Hardy-Sobolev inequalities with subcritical exponents}
\label{S:subcr}

Recall the Hardy inequality for fractional Laplacian:
\begin{equation}
\label{eq:Hardy}
\langle \Ds u,u\rangle\ge \mathcal{H}_s
\irn |x|^{-2s}|u|^2~\!dx
\end{equation}
for any $u\in \mathcal D^s(\R^n)$. The sharp constant in (\ref{eq:Hardy}) was found by Herbst in \cite{He}.

Next, take $\sigma\in(0,s)$. H\"older interpolation between the Sobolev and the Hardy inequalities
gives
\begin{equation}
\label{eq:Hardy-Sob}
\mathcal S_{s,\sigma}(\R^n)
:=\inf_{\scriptstyle u\in {\mathcal D}^s(\R^n)\atop\scriptstyle  u\ne 0}
\frac{\langle \Ds u,u\rangle}
{\||x|^{\sigma-s}u\|^2_{L^{\sigstar}(\R^n)}}>0.
\end{equation}
For $s=1$ and $n\ge 3$ the sharp value of $\mathcal S_{1,\sigma}(\R^n)$ was established in \cite{GY}.
It turns out that it is achieved by the function  
$$
U^{(\sigma)}(x)=\big(1+|x|^{\frac {2\sigma(n-2)}{n-2\sigma}}\big)^{\frac{2\sigma-n}{2\sigma}}.
$$
In the fractional case the next existence result holds.

\begin{Lemma}%[Subcritical case]
\label{L:spherical}
Let $s\in(0,1)$, $n>2s$, and $\sigma\in(0,s)$. Then the infimum $\mathcal S_{s,\sigma}(\R^n)$ is
achieved. 
\end{Lemma}

We skip the proof of Lemma \ref{L:spherical} 
because it can be obtained by adapting the argument we used in \cite[Theorem 1.1]{MNcn}. 
An alternative approach is to use the duality of (\ref{eq:Hardy-Sob}) and the weighted Hardy-Littlewood-Sobolev
inequality \cite{SW}. The attainability of the sharp constant in the last one was proved in \cite{Lb}.

Next we deal with Neumann Laplacians on half-spaces. By using the even extension
${\mathcal D}^s_{\rm R}(\R^n_+)\to \mathcal D^s(\R^n)$,
see (\ref{eq:even_ex}), and since trivially $\langle \DsNS u,u\rangle\ge \langle \DsNR u,u\rangle$, one plainly gets that 
the constants
$$
\mathcal S^{\rm R}_{s,\sigma}(\R^n_+)
:=\inf_{\scriptstyle u\in {\mathcal D}^s_{\rm R}(\R^n_+)\atop\scriptstyle  u\ne 0}
\frac{\langle \DsNR u,u\rangle}
{\||x|^{\sigma-s}u\|^2_{L^{\sigstar}(\R^n_+)}}~,\quad 
\mathcal S^{\rm Sr}_{s,\sigma}(\R^n_+)
:=\inf_{\scriptstyle u\in {\mathcal D}^s_{\rm Sr}(\R^n_+)\atop\scriptstyle  u\ne 0}
\frac{\langle \DsNS u,u\rangle}
{\||x|^{\sigma-s}u\|^2_{L^{\sigstar}(\R^n_+)}}
$$
are positive. The next existence results can be obtained by adapting the argument we used in the proofs of Theorems \ref{T:Restricted}, 
\ref{T:Semirestricted}. Notice that the assumption $\sigma<s$ implies $\sigstar<\sstar$ and 
the compactness of the embedding $\mathcal D^s(\R^n)\hookrightarrow L^\sigstar_{\rm loc}(\R^n)$.
This considerably simplifies the proof compared with Sections \ref{S:Restricted} and \ref{S:Semirestricted}. In particular, we do not need to
prove preliminary inequalities between sharp constants.

\begin{Theorem}%[Subcritical case]
\label{T:sphericalR_Sr}
Let  $s\in(0,1)$, $n>2s$, and $\sigma\in(0,s)$. 
\begin{itemize}
\item[$i)$] The infimum $\mathcal S^{\rm R}_{s,\sigma}(\R^n_+)$ is achieved in
$ {\mathcal D}^s_{\rm R}(\R^n_+)$;
\item[$ii)$] The infimum $\mathcal S^{\rm Sr}_{s,\sigma}(\R^n_+)$ is achieved in
$ {\mathcal D}^s_{\rm Sr}(\R^n_+)$ by a function $u$ such that $u=P_{\!s}u$.
\end{itemize}
\end{Theorem}

The following theorem is proved exactly as Theorem \ref{T:Spectral}.

\begin{Theorem}
\label{T:sphericalSp}
Let  $s\in(0,1)$, $n>2s$, and $\sigma\in(0,s)$. Then
$$
\mathcal S^{\rm Sp}_{s,\sigma}(\R^n_+):=\inf_{u\in {\mathcal D}^s_{\rm R}(\R^n_+)\atop\scriptstyle  u\ne 0}
\frac{\langle\DsNSp u,u\rangle}{\||x|^{\sigma-s}u\|^2_{L^{\sigstar}(\R^n_+)}}
=2^{-\frac{2\sigma}{n}}\mathcal S_{s,\sigma}(\R^n),
$$
and $\mathcal S^{\rm Sp}_{s,\sigma}(\R^n_+)$ is achieved. 
\end{Theorem}

\appendix

\section{The operator $P_{\!s}$}
\label{S:P}

We start with few general results of independent interest about the linear operator 
$$
({\mathcal P}_{\!s}u)(x)={\frac{C_{n,s}}{\gamma_s} |x_1|^{2s}\displaystyle\irnplus\frac{u(y)}{|x-y|^{n+2s}}~\!dy}~,\quad x\in\R^n_-~,
\quad u: \R^n_+\to\R.
$$
Let us define
$$
\mathcal B_s(\beta)=\frac{\Gamma\big(1-\beta\big)\Gamma\big(2s+\beta\big)}{\Gamma\big(2s\big)}~,\quad
-2s<\beta<1.
$$
Notice that
\begin{eqnarray*}
\frac{\gamma_s}{C_{n,s}}\mathcal B_s(\beta)=|x_1|^{2s+\beta}\irnplus\frac{y_1^{-\beta}dy}{|x-y|^{n+2s}}=
y_1^{2s+\beta}\irnminus\frac{|x_1|^{-\beta}dx}{|x-y|^{n+2s}}
\quad\forall x\in\R^n_-~, ~~\forall y\in\R^n_+.
\end{eqnarray*}
\begin{Lemma}
\label{L:convolution0}
Let $p\in(1,\infty)$, $t\in (-\frac{1}{p'},2s+\frac1p)$ and let $\alpha$ be an exponent satisfying
$$
-2s<\frac{\alpha}{p-1}<1~,\qquad 0<\alpha+tp<1+2s~\!.
$$
If $u\in L^p(\R^n_+;|x_1|^{-tp}dx)$, then ${\mathcal P}_{\!s}u\in L^p(\R^n_-;|x_1|^{-tp}dx)$ and
\begin{equation}
\label{eq:alpha0}
\irnminus|x_1|^{-tp} |{\mathcal P}_{\!s} u(x)|^{p} dx
\le \mathcal B_s\Big(\frac{\alpha}{p-1}\Big)^{p-1}\mathcal B_s(\alpha+tp-2s\big) 
~\!\irnplus |x_1|^{-tp}|u|^p~\!dx.
\end{equation}
\end{Lemma}

\proof
We use H\"older inequality and Fubini's theorem to estimate
\begin{eqnarray*}
\irnminus|x_1|^{-tp} |{\mathcal P}_{\!s}u(x)|^p dx\!\!\!&=&\!\!\!\left(\frac{C_{n,s}}{\gamma_s}\right)^{\!p}
\irnminus |x_1|^{p(2s-t)}
\Big[\irnplus\left(\frac{y_1^\alpha|u(y)|^p}{|x-y|^{n+2s}}\right)^\frac{1}{p}\!\!
\frac{y_1^{-\frac{\alpha}p}}{|x-y|^{\frac{n+2s}{p'}}}~\!dy\Big]^p dx\\
&\le&\left(\frac{C_{n,s}}{\gamma_s}\right)^{\!p} \irnminus |x_1|^{p(2s-t)}
\Big(\irnplus \frac{y_1^\alpha|u(y)|^pdy}{|x-y|^{n+2s}}\Big)~\!
\Big(\irnplus \frac{y_1^{-\frac{\alpha}{p-1}}dy}{|x-y|^{n+2s}}\Big)^{p-1}\!dx\\
\!\!\!&=&\!\!\! \frac{C_{n,s}}{\gamma_s}~\!\mathcal B_s\Big(\frac{\alpha}{p-1}\Big)^{p-1}
\irnminus |x_1|^{2s-tp-\alpha}
\Big(\irnplus \frac{y_1^\alpha|u(y)|^pdy}{|x-y|^{n+2s}}\Big)dx\\
\!\!\!&=&\!\!\! \frac{C_{n,s}}{\gamma_s}~\!\mathcal B_s\Big(\frac{\alpha}{p-1}\Big)^{p-1}
\irnplus y_1^\alpha|u(y)|^p~\!
\Big(\irnminus \frac{|x_1|^{2s-tp-\alpha}dx}{|x-y|^{n+2s}}\Big)dy\\
\!\!\!&=&\!\!\! \mathcal B_s\Big(\frac{\alpha}{p-1}\Big)^{p-1}\mathcal B_s(\alpha+tp-2s)~\!\irnplus y_1^{-tp}|u(y)|^p~\!dy.
\end{eqnarray*}
The proof is complete.
\QED

\begin{Corollary}
\label{C:convolution1}
The linear transform ${\mathcal P}_{\!s}$ is continuous 
$L^p(\R^n_+)\to L^p(\R^n_-)$
for any exponent 
$p\in (1,\infty]$. 
\end{Corollary}

\proof
If $p=\infty$ it trivially holds that ${\mathcal P}_{\!s}:L^\infty(\R^n_+)\to L^\infty(\R^n_-)$ is nonexpansive.
To handle the case $p\in(1,\infty)$ take $t=0$ in Lemma \ref{L:convolution0} and conclude.
\QED

\begin{Remark}
\label{R:P}
The assumption $p>1$ is  necessary. 
Let $E\subset \R^n_+$ be any bounded measurable set of positive measure. 
Since
$\chi_E\in L^p(\R^n_+)$ for any $p\in [1,\infty]$, then
${\mathcal P}_{\!s}\chi_E\in L^p(\R^n_-)$ for any $p\in (1,\infty]$ by Corollary \ref{C:convolution1}. 
Now, for $x\in\R^n_-$
and $R>0$ large enough we estimate
$$
({\mathcal P}_{\!s}\chi_E)(x)= c~\! |x_1|^{2s}\int\limits_E\frac{dy}{(y-x)^{n+2s}}
\ge c(E)~\! \frac{|x_1|^{2s}}{(R+|x|)^{n+2s}},
$$
that readily implies ${\mathcal P}_{\!s}\chi_E\notin L^1(\R^n_-)$.
\end{Remark}

\noindent
{\bf Proof of Lemma \ref{L:extension}.}
To check $i)$  note that the set 
$$K_u=\{ \omega\in \mathcal D^s_{\rm Sr}(\R^n_+)~|~\omega|_{\R^n_+}=u~~\text{on $\R^n_+$}~\}$$
is convex, closed and not empty.
Thus the minimization problem (\ref{eq:P_minimization}) has a unique solution $P_{\!s}u \in \mathcal D^s_{\rm Sr}(\R^n_+)$.
Further, we have that
\begin{equation}
\label{eq:Ort}
\langle\DsNS(P_{\!s}u),\f\rangle=0\quad\text{for any 
$\f \in \mathcal D^s_{\rm Sr}(\R^n_+)$ such that $\f\equiv 0$ on $\R^n_+$},
\end{equation}
because the polynomial 
$t\mapsto \mathcal E_s( {P_{\!s}} u+t\f;\R^{2n}\setminus(\R^n_-)^2)$ attains its minimum at $t=0$.
Take $\f\in\mathcal C^\infty_0(\R^n_-)$. Using (\ref{eq:Ort}), (\ref{eq:Sr_diff}) and
recalling that $\f(x)-\f(y)\equiv 0$ for $x,y\in \R^n_+$, we find
\begin{multline*}
0=C_{n,s}\int\limits_{\R^n_-}\f(x)\Big(\irnplus \frac{(P_{\!s}u)(x)-(P_{\!s}u)(y)}{|x-y|^{n+2s}}~\!\!dy\Big)dx
\\
=C_{n,s}\int\limits_{\R^n_-}\f(x)(P_{\!s}u)(x)\Big(\irnplus \frac{dy}{|x-y|^{n+2s}}\Big)dx-
C_{n,s}\int\limits_{\R^n_-}\f(x)\Big(\irnplus \frac{u(y)}{|x-y|^{n+2s}}~\!\!dy\Big)dx\\
=\int\limits_{\R^n_-}\f(x)\Big(\gamma_s|x_1|^{-2s} (P_{\!s}u)(x)-
C_{n,s}\irnplus \frac{u(y)}{|x-y|^{n+2s}}~\!\!dy\Big)dx~\!.
\end{multline*}
Since $\f$ was arbitrarily chosen, the identity (\ref{eq:P_explicit}) follows.

Now we prove $ii)$. The operator $u\mapsto {P_{\!s}} u$ is clearly linear. 
From $K_{\!{P_{\!s}} u}=K_u$ we infer that ${P_{\!s}}$ is projector. Since 
$$
Ker ({P_{\!s}})=\{\f \in \mathcal D^s_{\rm Sr}(\R^n_+)\ \ \big|\ \ \f\equiv 0\ \ {\rm on}\ \  \R^n_+\},
$$
we see from (\ref{eq:Ort}) that $Ker ({P_{\!s}})\perp Im ({P_{\!s}})$, thus ${P_{\!s}}$ is orthoprojector.

Finally, statement $iii)$  is the Pythagorean theorem for orthoprojectors.
\QED

\noindent
{\bf Proof of Lemma \ref{L:Pcompactness}.}
By Corollary \ref{C:convolution1} and thanks to Lemma \ref{L:Sr_space}, for any $u\in \mathcal D^s_{\rm Sr}(\R^n_+)$ we have
$$\|P_{\!s}u\|^2_{L^\sstar(\R^n)}\le c \|u\|^2_{L^\sstar(\R^n_+)}\le c ~\!\mathcal E_s(u;\R^{2n}\setminus(\R^n_-)^2).
$$
Thus $P_{\!s}:\mathcal D^s_{\rm Sr}(\R^n_+)\to L^\sstar(\R^n)$ is continuous. Since $P_{\!s}$ coincides with the identity on
$\mathcal R^s_{\rm Sr}(\R^n_+)$, the continuity of the embedding 
$\mathcal R^s_{\rm Sr}(\R^n_+)\hookrightarrow L^\sstar(\R^n)$ follows for free.

To prove $ii)$ take an exponent $p\in[1,\sstar)$ and a sequence 
$u_h\in \mathcal D^s_{\rm Sr}(\R^n_+)$ such that $u_h\to 0$ weakly in $\mathcal D^s_{\rm Sr}(\R^n_+)$. 
We have to show that $P_{\!s} u_h\to 0$ in $L^p(B_r)$ for any $r>0$.

For arbitrary $\rho>2r$ we write 
$$
 P_{\!s}u_h=  P_{\!s}(u_h\chi_{B_\rho})+ P_{\!s}(u_h\chi_{\R^n\setminus B_\rho})=:{\cal U}^0_{h}+{\cal U}^\infty_{h}.
$$
We estimate  ${\cal U}^\infty_h(x)$ for $x\in B_r^-$ as follows:
\begin{eqnarray*}
|{\cal U}^\infty_h(x)| &\le& c\, |x_1|^{2s}\!\int\limits_{\R^n_+\setminus B_\rho}\frac{|u_h(y)|\,dy}{|x-y|^{n+2s}}\\
&\le& c\, r^{2s} \|u_h\|_{L^\sstar(\R^n_+)} 
\Big(\int\limits_{\R^n\setminus B_\rho}\frac{dz}{|z|^{2n}}\Big)^{\frac{n+2s}{2n}}\!=
c\,\|u_h\|_{L^\sstar(\R^n_+)}\cdot\frac {r^{2s}}{\rho^{\frac{n+2s}{2}}}.
\end{eqnarray*}
We infer that
$$
\|{\cal U}_h^\infty\|_{L^p(B_r^-)}\le 
c\,\|u_h\|_{L^\sstar(\R^n_+)}\cdot \frac {r^{\frac{n}{p}+2s}}{\rho^{\frac{n+2s}{2}}}.
$$
Trivially 
$$
\|P_{\!s} u_h\|_{L^p(B_r)}\le \|u_h\|_{L^p(B_r^+)}+\|{\cal U}^0_h\|_{L^p(B_r^-)}+\|{\cal U}_h^\infty\|_{L^p(B_r^-)},
$$
Corollary \ref{C:convolution1} gives $\|{\cal U}^0_h\|_{L^p(B_r^-)}\le c(p)\,\|u_h\|_{L^p(B_\rho^+)}$, 
 so we arrive at
\begin{equation}
\label{eq:hot}
\|P_{\!s} u_h\|_{L^p(B_r)}\le  
c(p)\,\|u_h\|_{L^p(B_\rho^+)}+ c\,\|u_h\|_{L^\sstar(\R^n_+)} \cdot  \frac {r^{\frac{n}{p}+2s}}{\rho^{\frac{n+2s}{2}}}.
\end{equation}
Since $u_h\to 0$ weakly in $\mathcal D^s_{\rm Sr}(\R^n_+)$, we have that $u_h$ is bounded in $L^\sstar(\R^n_+)$.
Thus, given any $\eps>0$ we can find a large $\rho=\rho(\eps)>0$ such that the last term in (\ref{eq:hot}) is smaller than $\eps$. Hence
$$
\limsup_{h\to\infty}\|P_{\!s} u_h\|_{L^p(B_r)}\le  
c(p)\, \limsup_{h\to\infty}\|u_h\|_{L^p(B_\rho^+)}+ \eps=\eps,
$$
as by Lemmata \ref{L:Sr_space} and \ref{L:spaceR} 
we have $u_h\to 0$ in $L^p(B_\rho^+)$. Since $\eps>0$ was arbitrarily chosen,
we are done.
\QED

\noindent{\bf Proof of Lemma \ref{L:differential}.}
Let ${\mathcal F}$ be the Fourier transform in $\R^n$. It is well known that 
${\mathcal F}[\Ds\f](\xi)=|\xi|^{2s}{\mathcal F}[\f]$
for any $\f\in {\mathcal C}^\infty_0(\R^n)$. Since $u\in L^\sstar(\R^n)$
we can define the distribution  $\Ds u$ via
$$
\langle\Ds u,\f\rangle:=
\irn u(x)\Ds\f(x)~\!dx=
\irn|\xi|^{2s}\mathcal F[u]\overline{\mathcal F[\f]}~\!d\xi,\quad\f\in \mathcal C^\infty_0(\R^n) .
$$
Next, if  $\f\in  \mathcal C^\infty_0(\R^n_+)$
then $\f(x)-\f(y)=0$ for $x,y\in\R^n_-$, hence
$$
\langle\DsNS u,\f\rangle=
\frac{C_{n,s}}{2}\iirn\frac{(u(x)-u(y))(\f(x)-\f(y))}{|x-y|^{n+2s}}~\!\!dxdy=
\irn u(x)\Ds\f(x)~\!dx,
$$
that concludes the proof of $i)$.

Next, since $u=P_{\!s}u$ and $\mathcal C^\infty_0(\R^n_-)\subset \mathcal D^s_{\rm Sr}(\R^n_+)$,  then (\ref{eq:Ort}) gives $\DsNS u=0$ on $\R^n_-$ 
immediately. Using again $u=P_{\!s}u$ and the explicit expression for $P_{\!s}$ in (\ref{eq:P_explicit}),  for any $x\in \R^n_-$ we obtain
\begin{multline*}
0=\gamma_s|x_1|^{-2s}\big(u(x)-(P_{\!s}u)(x)\big)\\={C_{n,s}}\Big(\irnplus\frac{u(x)}{|x-y|^{n+2s}}~\!dy
-\irnplus\frac{u(y)}{|x-y|^{n+2s}}~\!dy\Big)
= \DsNR u(x),
\end{multline*}
and the lemma is completely proved.
\QED

\section{Limits}
\label{A:A}
The well known behaviors of $C_{n,s}$ as $s\to 0^+$ and $s\to 1^-$ follow from the 
 identity
$$
\frac{C_{n,s}}{s(1-s)}=
\frac{2^{2s}\Gamma\big(\frac{n}{2}+s\big)}{\pi^{\frac{n}{2}}\Gamma\big(2-s\big)}.
$$

Next, fix a  function $u\in H^1(\R^n)$. As in the proof of Lemma \ref{L:differential}
we denote by $\mathcal F$ the Fourier transform. Lebesgue's dominated convergence theorem and 
the classical identity
$$
\mathcal E_s(u;\R^n\!\times\R^n\!)=
\irn\left(|\xi|^{2}|\mathcal F[u]|^2\right)^s
|\mathcal F[u]|^{2(1-s)}~\!d\xi
$$
readily give
\begin{equation}
\label{eq:tutto}
\lim_{s\to 1^-}\mathcal E_s(u;\R^n\!\times\R^n\!)=\irn|\nabla u|^2~\!dx~,\quad
\lim_{s\to 0^+}\mathcal E_s(u;\R^n\!\times\R^n\!)=\irn| u|^2~\!dx.
\end{equation}
We are in position to compute the limits of the Neumann Restricted and Semirestricted quadratic forms on half-spaces, 
as $s(1-s)\to 0^+$.

\begin{Theorem}[Limits as $s\to 1^-$]
\label{P:s_to1}
If $u\in H^1(\R^n)$ then
$$
\lim_{s\to 1^-}\mathcal E_s(u;\R^n_+\!\times\R^n_+\!)=\lim_{s\to 1^-}\mathcal E_s(u;\R^{2n}\setminus(\R^n_-)^2)= \irnplus |\nabla u|^2~\!dx.
$$
\end{Theorem}

\proof 
The conclusion for the Restricted quadratic form should be known, at least for bounded domains.
We cite for instance \cite{BBM} for related results.
We furnish here a complete proof for the convenience of the reader.

We denote by $c_n$ any constant possibly depending on the dimension $n$
but not on $s$; in particular, we have $C_{n,s}\le c_n~\!(1-s)$ for $s\in(0,1)$.
We start by proving that
\begin{equation}
\label{eq:claim1}
\mathcal E_s(u;\R^n_+\!\times\R^n_-\!)=o(1)\quad\text{as $s\to 1$.}
\end{equation}
{Note that the proof in \cite[Subsection 5.1]{DRoV}
of a similar result on bounded domains $\Omega\subset\R^n$ contains a defect, precisely in the proof of 
formula (5.5). However, the statement of \cite[Proposition 5.1]{DRoV} is correct.}

We introduce the notation $\Pi_\delta=\{x\in\R^n\,\big|\,|x_1|<\delta\}$
and estimate
$$
\frac{2\mathcal E_s(u;\R^n_+\!\times \R^n_-)}{C_{n,s}}\le 
\int\limits_{\R^n_+}dx\!\!\!\int\limits_{~\R^n_-\setminus B_\delta(x)}\!\!\!\frac{(u(x)-u(y))^2}{|x-y|^{n+2s}}\,dy+
\int\limits_{\Pi_\delta}dx\int\limits_{B_\delta(x)}\!\frac{(u(x)-u(y))^2}{|x-y|^{n+2s}}\,dy=: \mathcal{I}_1+\mathcal{I}_2.
$$
We have
$$
\mathcal{I}_1\le 2\iint\limits_{\{|x-y|>\delta\}}\frac{u(x)^2+u(y)^2}{|x-y|^{n+2s}}\,dxdy=
4\int\limits_{\R^n}u(x)^2\,dx\int\limits_{\{|z|>\delta\}}\frac{dz}{|z|^{n+2s}}\le 
c_n\,\|u\|_{L^2(\R^n)}^2\,
{\delta^{-2s}}.
$$
To handle $\mathcal{I}_2$ we estimate
$$
(u(y)-u(x))^2\!=\!\Big(\int\limits_0^1\nabla u(x+\tau(y-x))\cdot (y-x)\,d\tau\Big)^{\!2}\!\!\le
|x-y|^2\!\int\limits_0^1|\nabla u(x+\tau(y-x))|^2\,d\tau,
$$
so that 
\begin{eqnarray*}
\mathcal{I}_2 &\le& 
\int\limits_0^1d\tau\!\int\limits_{\Pi_\delta}dx\!\int\limits_{B_\delta(x)} \frac{|\nabla u(x+\tau(y-x))|^2}{|x-y|^{n+2s-2}}\,dy
=\int\limits_0^1d\tau\!\int\limits_{B_\delta}\frac{dz}{|z|^{n+2s-2}}
\int\limits_{\Pi_\delta} {|\nabla u(x+\tau z)|^2}\,dx\\
&\le&\int\limits_{B_\delta}\frac{dz}{|z|^{n+2s-2}}
\int\limits_{\Pi_{2\delta}} {|\nabla u(x)|^2}\,dx
= c_n\|\nabla u\|_{L^2(\Pi_{2\delta})}^2\,\frac {\delta^{2(1-s)}}{2(1-s)}.
\end{eqnarray*}
Thus, for $s$ close to $1$ and small $\delta$, we obtain 
$$
{\mathcal E_s(u;\R^n_+\!\times \R^n_-)}\le 
c_n (1-s)\,{\|u\|_{L^2(\R^n)}^2}\,\delta^{-2s}+c_n\|\nabla u\|_{L^2(\Pi_{2\delta})}^2.
$$
Formula  (\ref{eq:claim1}) readily follows, because 
for any 
$\eps>0$ we can find  $\delta=\delta(\eps)>0$ 
such that 
$c_n\|\nabla u\|_{L^2(\Pi_{2\delta})}^2<\eps$, and thus
$$
\limsup_{s\to 1^+} {\mathcal E_s(u;\R^n_+\!\times \R^n_-)}\le
\limsup_{s\to 1^+} c_n(1-s)\,{\|u\|_{L^2(\R^n)}^2}\,\delta^{-2s}+\eps
=
 \eps.
$$

From (\ref{eq:claim1}) we first infer that 
$\mathcal E_s(u;\R^n\!\times\R^n\!)=\mathcal E_s(u;\R^n_+\!\times\R^n_+\!)+
\mathcal E_s(u;\R^n_-\!\times\R^n_-\!)+o(1)$.
Further, by replacing $u$ by $\hat u$, 
that is the symmetric extension of $u|_{\R^n_+}$, and using (\ref{eq:tutto}) we obtain  
$$
2\irnplus|\nabla u|^2~\!dx=
\irn|\nabla \hat u|^2~\!dx=
\mathcal E_s(\hat u;\R^n\!\times\R^n\!)+o(1)=2\mathcal E_s(u;\R^n_+\!\times\R^n_+\!)+o(1).
$$
Thus $\displaystyle{\irnplus|\nabla u|^2~\!dx=\mathcal E_s(u;\R^n_+\!\times\R^n_+\!)+o(1)}$.
The conclusion for $\mathcal E_s(u;\R^{2n}\setminus(\R^n_-)^2)$ readily follows from 
(\ref{eq:all}) and (\ref{eq:claim1}).
\QED

Now we study the limits as $s\to 0$. It is convenient to discuss separately the
Restricted and the Semirestricted cases.

\begin{Theorem}[Limit as $s\to 0^+$, Restricted Laplacian]
\label{P:restricted0}
If $u\in H^1(\R^n)$ then
$$
\lim_{s\to 0^+}\mathcal E_s(u;\R^n_+\!\times\R^n_+\!)=\frac12\irnplus| u|^2~\!dx.
$$
\end{Theorem}

\proof
For $s\in(0,\frac12)$ we have $\chi_{\R^n_+}u\in H^s(\R^n)$, see
\cite[Sec. 2.10.2]{Tr}. Thus, via (\ref{eq:N_sum}) we can compute 
\begin{eqnarray*}
\mathcal E_s(u;\R^n_+\!\times\R^n_+)=
\mathcal E_s(\chi_{\R^n_+}u;\R^n_+\!\times\R^n_+)=
\mathcal E_s(\chi_{\R^n_+}u;\R^n\!\times\R^n)-\gamma_s\irnplus x_1^{-2s}|u|^2 dx.
\end{eqnarray*}
By (\ref{eq:tutto}) we have
$\displaystyle{\lim_{s\to 0^+}\mathcal E_s(\chi_{\R^n_+}u;\R^n\!\times\R^n)=
\irn|\chi_{\R^n_+}u|^2dx=\irnplus|u|^2dx}$. Next, from  (\ref{eq:gamma}) we see that $\displaystyle\lim_{s\to 0^+}\gamma_s=\frac12$, 
and the conclusion readily follows. 
\QED

\begin{Theorem}[Limit as $s\to 0^+$, Semirestricted Laplacian]
If $u\in H^1(\R^n)$ then
\begin{eqnarray}
\label{eq:Zu}
\lim_{s\to 0^+}\mathcal E_s(u;\R^{2n}\setminus(\R^n_-)^2)~&=&\irnplus | u|^2~\!dx+\frac12\irnminus | u|^2~\!dx,\\
\label{eq:ZPu}
\lim_{s\to 0^+}\!\mathcal E_s(P_{\!s}u;\R^{2n}\setminus(\R^n_-)^2)&=& \irnplus | u|^2~\!dx.
\end{eqnarray}
\end{Theorem}

\proof
Identity (\ref{eq:Zu}) readily  follows from (\ref{eq:few}) and Theorem \ref{P:restricted0}, since
$$
\mathcal E_s(u;\R^{2n}\setminus(\R^n_-)^2)=\mathcal E_s(u;\R^n\times\R^n)-
\mathcal E_s(u;\R^n_-\times\R^n_-)=
\irn  | u|^2~\!dx-\frac12\irnminus | u|^2~\!dx+o(1).
$$

To prove (\ref{eq:ZPu}) we first notice that $u\in L^2(\R^n_+; x_1^{-2s}dx)$ for $
0\le s<\frac12$. 
Choosing $p=2$, $t=s$ and $\alpha=1-\sqrt{s}$ in (\ref{eq:alpha0}), we get
that $P_{\!s}u\in L^2(\R^n_-; |x_1|^{-2s}dx)$ and 
$$
\irnminus|x_1|^{-2s} |P_{\!s} u|^{2} dx
\le \left(\frac{\Gamma\big(\sqrt{s}\big)\Gamma\big(1-\sqrt{s}+2s\big)}{\Gamma\big(2s\big)}\right)^2~\!\irnplus |x_1|^{-2s}|u|^2~\!dx.
$$
Since $\Gamma\big(1-\sqrt{s}+2s\big)=1+o(1)$, and
 $\displaystyle{\frac{\Gamma\big(\sqrt{s}\big)}{\Gamma\big(2s\big)}=O(\sqrt{s})}$, we
 infer that
\begin{equation}
\label{eq:ex_lemma}
|x_1|^{-s}P_{\!s} u\to 0\quad\text{in $L^2(\R^n_-)$. }
\end{equation}
Using also (\ref{eq:orthogonal}) we
obtain
\begin{eqnarray*}
\mathcal E_s(P_{\!s}u;\R^{2n}\setminus(\R^n_-)^2)&=& 
\mathcal E_s(u;\R^{2n}\setminus(\R^n_-)^2)-\mathcal E_s(u-P_{\!s}u;\R^{2n}\setminus(\R^n_-)^2)\\
&=&\mathcal E_s(u;\R^{2n}\setminus(\R^n_-)^2)-
\gamma_s\irnminus|x_1|^{-2s}|u-P_{\!s}u|^2~\!dx.
\end{eqnarray*}
The conclusion follows, thanks to  (\ref{eq:Zu}) and (\ref{eq:ex_lemma}).
 \QED

\footnotesize
%\section{references}
\label{References}


\begin{thebibliography}{XX}
\footnotesize

\bibitem{Au} 
T. Aubin, Probl\`emes isop\'erim\'etriques et espaces de Sobolev, J. Differential Geometry {\bf 11} (1976), no.~4, 573--598. 

\bibitem{AM1}
Adimurthi\ and\ G. Mancini, The Neumann problem for elliptic equations with critical nonlinearity, 
in {\it Nonlinear analysis}, 9--25, Sc. Norm. Super. di Pisa Quaderni, Scuola Norm. Sup., Pisa (1991).

\bibitem{AMY}
Adimurthi, G. Mancini and\ S.L. Yadava, The role of the mean curvature in semilinear Neumann problem involving critical exponent, Comm. Partial Differential Equations {\bf 20} (1995), no.~3-4, 591--631.


\bibitem{BBM}
J. Bourgain, H. Brezis\ and\ P. Mironescu, Another look at Sobolev spaces, in {\it Optimal control and partial differential equations}, 439--455, IOS, Amsterdam.

\bibitem{CoTa}
A. Cotsiolis\ and\ N. K. Tavoularis, Best constants for Sobolev inequalities for higher order 
fractional derivatives, J. Math. Anal. Appl. {\bf 295} (2004), no.~1, 225--236.

\bibitem{delMM}
M. del Pino, F. Mahmoudi\ and\ M. Musso, Bubbling on boundary submanifolds for the Lin-Ni-Takagi problem at higher critical exponents, J. Eur. Math. Soc. (JEMS) {\bf 16} (2014), no.~8, 1687--1748.

\bibitem{DRoV}
S. Dipierro, X. Ros-Oton\ and\ E. Valdinoci, Nonlocal problems with Neumann boundary conditions, Rev. Mat. Iberoam. {\bf 33} (2017), no.~2, 377--416.

\bibitem{DZ}
J. Dou\ and\ M. Zhu, Sharp Hardy-Littlewood-Sobolev inequality on the upper half space, Int. Math. Res. Not. {\bf 2015}, no.~3, 651--687.

\bibitem{DLV}
B. Dyda, J, Lehrb\"ack, A. V. V\"ah\"akangas,
Fractional Hardy-Sobolev type inequalities for half spaces and John domains,
preprint ArXiv:1709.03296 (2017).

\bibitem{Frank}
R.L. Frank, T. Jin\ and\ J. Xiong,
Minimizers for the fractional Sobolev inequality on domains,
preprint ArXiv 1707.00131 (2017).

\bibitem{GM}
M. Gazzini\ and\ R. Musina, 
On a Sobolev-type inequality related to the weighted $p$-Laplace operator, J. Math. Anal. Appl. {\bf 352} (2009), no.~1, 99--111.

\bibitem{GY}
N. Ghoussoub\ and\ C. Yuan, Multiple solutions for quasi-linear PDEs involving the critical Sobolev and Hardy exponents, 
Trans. Amer. Math. Soc. {\bf 352} (2000), no.~12, 5703--5743.

\bibitem{G1}
Q.-Y. Guan, Integration by parts formula for regional fractional Laplacian, Comm. Math. Phys. {\bf 266} (2006), no.~2, 289--329. 

\bibitem{G2}
Q. Y. Guan and Z. M. Ma, Boundary problems for fractional Laplacians, Stoch. Dyn., {\bf 5}
(2005), 385--424.

\bibitem{He}
I. W. Herbst, Spectral theory of the operator $(p\sp{2}+m\sp{2})\sp{1/2}-Ze\sp{2}/r$, Comm. Math. 
Phys. {\bf 53} (1977), no.~3, 285--294.

\bibitem{KN}
N. Kuznetsov\ and\ A. Nazarov, Sharp constants in the Poincar\'e, Steklov and related inequalities (a survey), Mathematika {\bf 61} (2015), no.~2, 328--344.

\bibitem{Lb} 
E. H. Lieb, Sharp constants in the Hardy--Littlewood--Sobolev and related inequalities, Ann. Math. {\bf 118} (1983), no.~2, 349--374. 

\bibitem{MNmp}
R. Musina\ and\ A. I. Nazarov, Some strong maximum principles for fractional Laplacians, 
Proc. Roy. Soc. Edinburgh Sect. A,
to appear. %preprint arXiv:1612.01043 (2016).

\bibitem{MNcn}
R. Musina\ and\ A. I. Nazarov, Fractional Hardy-Sobolev inequalities on half spaces, preprint 
arXiv:1707.02710 (2017).

\bibitem{Naz-surv}
A.I. Nazarov, Dirichlet and Neumann problems to critical Emden--Fowler type equations, J. Global Optim. {\bf40} (2008), 289--303.

\bibitem{Rosen}
G. Rosen, Minimum value for $c$ in the Sobolev inequality $\|\phi^3 \| \le c \|\nabla \phi \|^3$, SIAM J. Appl. Math. {\bf 21} (1971), 30--32.

\bibitem{Ro}
X. Ros-Oton, Nonlocal elliptic equations in bounded domains: a survey, Publ. Mat. {\bf 60} (2016), no.~1, 3--26.

\bibitem{SW}
E. M. Stein\ and\ G. Weiss, Fractional Integrals on $n$-dimensional Euclidean Space, Journ. Math. Mech. {\bf 7} (1958),  no.~4, 503--514. 

\bibitem{Ta}
G. Talenti, Best constant in Sobolev inequality, Ann. Mat. Pura Appl. (4) {\bf 110} (1976), 353--372.

\bibitem{Tr}
H. Triebel, Interpolation theory, function spaces, differential operators, Deutscher Verlag 
Wissensch., Berlin, 1978. 

\bibitem{WWY}
L. Wang, J. Wei\ and\ S. Yan, A Neumann problem with critical exponent in nonconvex domains and Lin-Ni's conjecture, 
Trans. Amer. Math. Soc. {\bf 362} (2010), no.~9, 4581--4615. 

\bibitem{W1}
M. Warma, The $p$-Laplace operator with the nonlocal Robin boundary conditions on arbitrary open sets, 
Ann. Mat. Pura Appl. (4) {\bf 193} (2014), no.~1, 203--235.

\bibitem{W2}
M. Warma, The fractional relative capacity and the fractional Laplacian 
with Neumann and Robin boundary conditions on open sets, Potential Anal. {\bf 42} (2015), no.~2, 499--547. 

\bibitem{W3}
M. Warma, A fractional Dirichlet-to-Neumann operator on bounded Lipschitz domains, 
Commun. Pure Appl. Anal. {\bf 14} (2015), no.~5, 2043--2067.

\end{thebibliography}
\end{document}